\documentclass{amsart}
\usepackage{amsmath, amsthm, amsfonts}
\usepackage{pgf,tikz}
\usepackage{mathrsfs}
\usepackage{hyperref}
\usetikzlibrary{arrows}

\newtheorem{conj}{Conjecture}
\newtheorem{prop}{Proposition}
\newtheorem{cor}{Corollary}
\newtheorem{thm}{Theorem}
\newtheorem{lem}{Lemma}

\theoremstyle{remark}

\theoremstyle{definition}
\newtheorem{defi}{Definition}

\newcommand{\Hlam}{H_{\lambda}}
\newcommand{\Alam}{P_{\lambda}}
\newcommand{\Blam}{Q_{\lambda}}
\newcommand{\Clam}{R_{\lambda}}
\newcommand{\RR}{\mathbb{R}}
\newcommand{\NN}{\mathbb{N}}
\newcommand{\CC}{\mathbb{C}}

\begin{document}

\title{A family of projective heat maps}
\author{Quang-Nhat Le}
\thanks{This project has received partial funding from the European Research Council (ERC) under the European Union's Horizon 2020 research and innovation programme (grant agreement No [ERC StG 716424 - CASe]).}

\address{Department of Mathematics, Brown University, 
	Box 1917, 151 Thayer Street, Providence, RI 02912 and Einstein Institute of Mathematics, Hebrew University of Jerusalem, Jerusalem 9190401, Israel}
\email{qnhatle@math.brown.edu, qnhatle@math.huji.ac.il}

\maketitle

\begin{abstract}
	The pentagram map was invented by Richard Schwartz in his search for a projective-geometric analogue of the midpoint map. 
	It turns out that the dynamical behavior of the pentagram map is totally different from that of the midpoint map. 
	Recently, Schwartz has constructed a related map, the projective heat map, which empirically exhibits similar dynamics 
	as the midpoint map. In this paper, we will demonstrate that there is a one-parameter family of maps which behaves a lot like Schwartz'
	projective heat map.
\end{abstract}

\section{Introduction}
The midpoint map is perhaps the simplest polygon iteration. Starting with an $n$-gon $P_1$, we create a new $n$-gon
$P_2$  whose vertices are midpoints of the edges of $P_1$. For almost every choice of $P_1$, if we iterate this process,
the obtained sequence of polygons $\{P_k\}$ will converge to a point. Furthermore, we can rescale the polygons $P_k$ so
that the new sequence converges to an affinely regular polygon.

The proof of these claims is elegant. The midpoint map can be thought of as a $\CC$-linear map whose
eigenvectors  are $v_j = (\omega^j, ..., \omega^{nj})$. The largest eigenvalue (in absolute value) is $1$ with
eigenvector $v_0$, which justifies the convergence of $\{P_k\}$ to a point. The eigenvectors corresponding the second
largest eigenvalue are $v_1$ and $v_{n-1}$. Because any complex linear combination of $v_1$ and $v_{n-1}$ is affinely regular, this implies the rescaled sequence $c_k P_k$ (with $c_k \in \mathbb{C}$) 
converges to a regular polygon.

The above arguments can be cast in the theory of discrete Fourier transform; see \cite{Te}. The midpoint map also 
has connection with the heat equation (\cite{Te}) and outer billiards (\cite{Ta} and \cite{Tr}).

In his first attempt to find a projectively natural analogue of the midpoint map, Richard Schwartz has invented the 
pentagram map, which is defined as follows. Starting with a \textit{generic} $n$-gon $P$ (i.e. none of the triples of
vertices lies on a line), we draw its shortest diagonals. Then the intersections of these diagonals form a new $n$-gon $Q$.

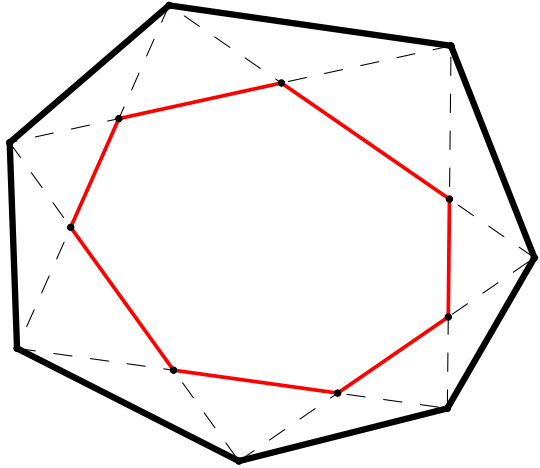
\begin{figure}[!ht]
	\centering
	\begin{tikzpicture}[line cap=round,line join=round,>=triangle 45,x=1.0cm,y=1.0cm, scale=0.8]
	\draw [line width=2.4pt] (-2.6021129085678947,0.5989637701741396)-- (0.046226113444625745,2.880609696831082);
	\draw [line width=2.4pt] (0.046226113444625745,2.880609696831082)-- (4.7317489985437,2.2083390220125185);
	\draw [line width=2.4pt] (4.7317489985437,2.2083390220125185)-- (6.117034025442558,-1.315989061127223);
	\draw [line width=2.4pt] (6.117034025442558,-1.315989061127223)-- (4.670633482651104,-3.8217252127236865);
	\draw [line width=2.4pt] (4.670633482651104,-3.8217252127236865)-- (1.2074209154039615,-4.697714273850905);
	\draw [line width=2.4pt] (1.2074209154039615,-4.697714273850905)-- (-2.4798818767827013,-2.823505119811274);
	\draw [line width=2.4pt] (-2.6021129085678947,0.5989637701741396)-- (-2.4798818767827013,-2.823505119811274);
	\draw [dash pattern=on 7pt off 7pt] (-2.4798818767827013,-2.823505119811274)-- (4.670633482651104,-3.8217252127236865);
	\draw [dash pattern=on 7pt off 7pt] (4.670633482651104,-3.8217252127236865)-- (4.7317489985437,2.2083390220125185);
	\draw [dash pattern=on 7pt off 7pt] (4.7317489985437,2.2083390220125185)-- (-2.6021129085678947,0.5989637701741396);
	\draw [dash pattern=on 7pt off 7pt] (-2.6021129085678947,0.5989637701741396)-- (1.2074209154039615,-4.697714273850905);
	\draw [dash pattern=on 7pt off 7pt] (1.2074209154039615,-4.697714273850905)-- (6.117034025442558,-1.315989061127223);
	\draw [dash pattern=on 7pt off 7pt] (6.117034025442558,-1.315989061127223)-- (0.046226113444625745,2.880609696831082);
	\draw [dash pattern=on 7pt off 7pt] (0.046226113444625745,2.880609696831082)-- (-2.4798818767827013,-2.823505119811274);
	\draw [red, line width=1.4pt] (-1.588398659806022,-0.8104785008316736)-- (-0.7879080466618661,0.9970809482035184);
	\draw [red, line width=1.4pt] (-0.7879080466618661,0.9970809482035184)-- (1.9134118185022277,1.589870585281195);
	\draw [red, line width=1.4pt] (1.9134118185022277,1.589870585281195)-- (4.70591598614076,-0.3405182017441019);
	\draw [red, line width=1.4pt] (4.70591598614076,-0.3405182017441019)-- (4.686039617314619,-2.3016532592568413);
	\draw [red, line width=1.4pt] (4.686039617314619,-2.3016532592568413)-- (2.848489311785467,-3.56735180995327);
	\draw [red, line width=1.4pt] (2.848489311785467,-3.56735180995327)-- (0.12052723160439807,-3.1865251947712965);
	\draw [red, line width=1.4pt] (0.12052723160439807,-3.1865251947712965)-- (-1.588398659806022,-0.8104785008316736);
	\begin{scriptsize}
	\draw [fill=black] (-2.6021129085678947,0.5989637701741396) circle (1.5pt);
	\draw [fill=black] (0.046226113444625745,2.880609696831082) circle (1.5pt);
	\draw [fill=black] (4.7317489985437,2.2083390220125185) circle (1.5pt);
	\draw [fill=black] (6.117034025442558,-1.315989061127223) circle (1.5pt);
	\draw [fill=black] (4.670633482651104,-3.8217252127236865) circle (1.5pt);
	\draw [fill=black] (1.2074209154039615,-4.697714273850905) circle (1.5pt);
	\draw [fill=black] (-2.4798818767827013,-2.823505119811274) circle (1.5pt);
	\draw [fill=black] (-1.588398659806022,-0.8104785008316736) circle (1.5pt);
	\draw [fill=black] (-0.7879080466618661,0.9970809482035184) circle (1.5pt);
	\draw [fill=black] (1.9134118185022277,1.589870585281195) circle (1.5pt);
	\draw [fill=black] (4.70591598614076,-0.3405182017441019) circle (1.5pt);
	\draw [fill=black] (4.686039617314619,-2.3016532592568413) circle (1.5pt);
	\draw [fill=black] (2.848489311785467,-3.56735180995327) circle (1.5pt);
	\draw [fill=black] (0.12052723160439807,-3.1865251947712965) circle (1.5pt);
	\end{scriptsize}
	\end{tikzpicture}
	\caption{The pentagram map.} \label{fig:penta}
\end{figure}
	
However, the pentagram map is actually a discrete integrable system (\cite{OST1} and \cite{OST2}). Its typical orbit closure is a torus, 
which can be given a flat structure so that the pentagram map acts as a translation on the orbit closure. The pentagram 
map has been studied and generalized by various authors (\cite{GP}, \cite{GSTV}, \cite{KS} and \cite{M-B}). It also has a surprising intimate link to cluster
algebras (\cite{Gl} and \cite{GP}).

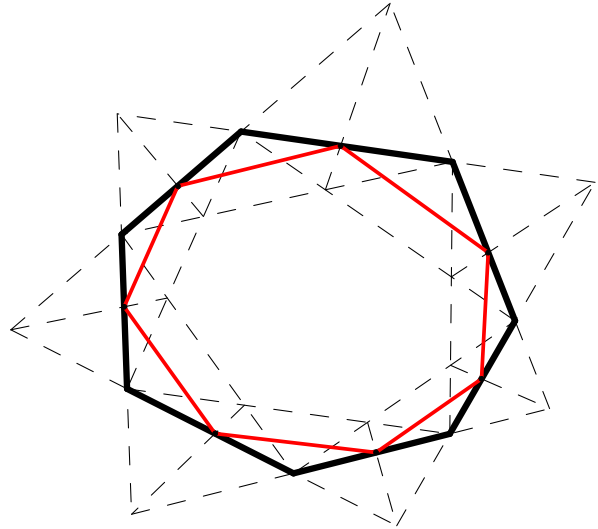
\begin{figure}[!ht]
	\centering
	\begin{tikzpicture}[line cap=round,line join=round,>=triangle 45,x=1.0cm,y=1.0cm, scale=0.6]
	\draw [line width=2.4pt] (-2.6021129085678947,0.5989637701741396)-- (0.046226113444625745,2.880609696831082);
	\draw [line width=2.4pt] (0.046226113444625745,2.880609696831082)-- (4.7317489985437,2.2083390220125185);
	\draw [line width=2.4pt] (4.7317489985437,2.2083390220125185)-- (6.117034025442558,-1.315989061127223);
	\draw [line width=2.4pt] (6.117034025442558,-1.315989061127223)-- (4.670633482651104,-3.8217252127236865);
	\draw [line width=2.4pt] (4.670633482651104,-3.8217252127236865)-- (1.2074209154039615,-4.697714273850905);
	\draw [line width=2.4pt] (1.2074209154039615,-4.697714273850905)-- (-2.4798818767827013,-2.823505119811274);
	\draw [line width=2.4pt] (-2.6021129085678947,0.5989637701741396)-- (-2.4798818767827013,-2.823505119811274);
	\draw [dash pattern=on 7pt off 7pt] (-2.4798818767827013,-2.823505119811274)-- (4.670633482651104,-3.8217252127236865);
	\draw [dash pattern=on 7pt off 7pt] (4.670633482651104,-3.8217252127236865)-- (4.7317489985437,2.2083390220125185);
	\draw [dash pattern=on 7pt off 7pt] (4.7317489985437,2.2083390220125185)-- (-2.6021129085678947,0.5989637701741396);
	\draw [dash pattern=on 7pt off 7pt] (-2.6021129085678947,0.5989637701741396)-- (1.2074209154039615,-4.697714273850905);
	\draw [dash pattern=on 7pt off 7pt] (1.2074209154039615,-4.697714273850905)-- (6.117034025442558,-1.315989061127223);
	\draw [dash pattern=on 7pt off 7pt] (6.117034025442558,-1.315989061127223)-- (0.046226113444625745,2.880609696831082);
	\draw [dash pattern=on 7pt off 7pt] (0.046226113444625745,2.880609696831082)-- (-2.4798818767827013,-2.823505119811274);
	\draw [dash pattern=on 7pt off 7pt] (-5.055228130744793,-1.5144893443167295)-- (3.4940736662314396,-5.859990810183104);
	\draw [dash pattern=on 7pt off 7pt] (3.4940736662314396,-5.859990810183104)-- (7.889846181602513,1.7552207305301668);
	\draw [dash pattern=on 7pt off 7pt] (7.889846181602513,1.7552207305301668)-- (-2.697660551904481,3.2742977835985627);
	\draw [dash pattern=on 7pt off 7pt] (6.882081180923875,-3.262359030219397)-- (3.3490407324734592,5.726111522455926);
	\draw [dash pattern=on 7pt off 7pt] (-5.055228130744793,-1.5144893443167295)-- (3.3490407324734592,5.726111522455926);
	\draw [dash pattern=on 7pt off 7pt] (-5.055228130744793,-1.5144893443167295)-- (-1.588398659806022,-0.8104785008316736);
	\draw [dash pattern=on 7pt off 7pt] (-0.7879080466618661,0.9970809482035184)-- (-2.697660551904481,3.2742977835985627);
	\draw [dash pattern=on 7pt off 7pt] (1.9134118185022277,1.589870585281195)-- (3.3490407324734592,5.726111522455926);
	\draw [dash pattern=on 7pt off 7pt] (4.70591598614076,-0.3405182017441019)-- (7.889846181602513,1.7552207305301668);
	\draw [dash pattern=on 7pt off 7pt] (4.686039617314619,-2.3016532592568413)-- (6.882081180923875,-3.262359030219397);
	\draw [dash pattern=on 7pt off 7pt] (2.848489311785467,-3.56735180995327)-- (3.4940736662314396,-5.859990810183104);
	\draw [dash pattern=on 7pt off 7pt] (0.12052723160439807,-3.1865251947712965)-- (-2.380533641044462,-5.605255720481978);
	\draw [dash pattern=on 7pt off 7pt] (-2.697660551904481,3.2742977835985627)-- (-2.380533641044462,-5.605255720481978);
	\draw [dash pattern=on 7pt off 7pt] (-2.380533641044462,-5.605255720481978)-- (6.882081180923875,-3.262359030219397);
	\draw [red, line width=1.4pt] (-2.5448390828690397,-1.0047033493937998)-- (-1.3550534648790113,1.6733534447368705);
	\draw [red, line width=1.4pt] (-1.3550534648790113,1.6733534447368705)-- (2.2515834829052217,2.5641888568649964);
	\draw [red, line width=1.4pt] (2.2515834829052217,2.5641888568649964)-- (5.522374634527068,0.19689438929012487);
	\draw [red, line width=1.4pt] (5.522374634527068,0.19689438929012487)-- (5.374276812655812,-2.602737472011303);
	\draw [red, line width=1.4pt] (5.374276812655812,-2.602737472011303)-- (3.0365099599727006,-4.235062339048223);
	\draw [red, line width=1.4pt] (3.0365099599727006,-4.235062339048223)-- (-0.5292998632700097,-3.8149611708784996);
	\draw [red, line width=1.4pt] (-0.5292998632700097,-3.8149611708784996)-- (-2.5448390828690397,-1.0047033493937998);
	\begin{scriptsize}
	\draw [fill=black] (-2.6021129085678947,0.5989637701741396) circle (1.5pt);
	\draw [fill=black] (0.046226113444625745,2.880609696831082) circle (1.5pt);
	\draw [fill=black] (4.7317489985437,2.2083390220125185) circle (1.5pt);
	\draw [fill=black] (6.117034025442558,-1.315989061127223) circle (1.5pt);
	\draw [fill=black] (4.670633482651104,-3.8217252127236865) circle (1.5pt);
	\draw [fill=black] (1.2074209154039615,-4.697714273850905) circle (1.5pt);
	\draw [fill=black] (-2.4798818767827013,-2.823505119811274) circle (1.5pt);
	\draw [fill=black] (-2.5448390828690397,-1.0047033493937998) circle (1.5pt);
	\draw [fill=black] (-1.3550534648790113,1.6733534447368705) circle (1.5pt);
	\draw [fill=black] (2.2515834829052217,2.5641888568649964) circle (1.5pt);
	\draw [fill=black] (5.522374634527068,0.19689438929012487) circle (1.5pt);
	\draw [fill=black] (-0.5292998632700097,-3.8149611708784996) circle (1.5pt);
	\draw [fill=black] (3.0365099599727006,-4.235062339048223) circle (1.5pt);
	\draw [fill=black] (5.374276812655812,-2.602737472011303) circle (1.5pt);
	\end{scriptsize}
	\end{tikzpicture}
	\caption{Schwartz' projective heat map.} \label{fig:heat}
\end{figure}

Schwartz' second attempt in \cite{S1} gave rise to the projective heat map $H$, whose construction is a bit more complicated than 
that of the pentagram map, as shown in Figure \ref{fig:heat}. In words, we connect the points obtained from the pentagram map and its inverse and then intersect with edges of the original polygon $P$. The new polygon $Q$ is inscribed in $P$.

Computer simulation suggests that this map behaves much like the midpoint map: for almost every starting (generic)
$n$-gon,  the iteration of $H$ will converge, up to projective transformations, to a projectively 
regular $n$-gon; see Conjecture 1.1 in \cite{S1}. Unlike the midpoint map, the projective heat map is non-linear, which 
makes proving the convergence claim significantly harder. Only the case of pentagons is resolved so far and its proof 
relies heavily on computer calculation. Even in this single case, the set of pentagons which do not exhibit the above convergence 
behavior possesses beautiful symmetries and structures as shown in \cite{S1}.

We can further generalize the projective heat map as follows. First we group the points obtained from the projective 
heat map, the pentagram map and its inverse into collinear triples. We can parametrize the line through each triple by 
mapping it to the projective line $\RR P^1$ so that the points corresponding to the pentagram map, the projective heat
map and the inverse pentagram map are mapped to the points $0$, $1$ and $\infty$, respectively. Next, pick a parameter 
$\lambda \in {\RR} \cup \infty = {\RR P}^1$. This determines, for each triple, a fourth point lying on the same
line.  We call this construction the parametrized heat map $H_{\lambda}$. Now we have a $1$-parameter family of maps 
interpolating the projective heat map $H_1$, the pentagram map $H_0$ and its inverse $H_{\infty}$.

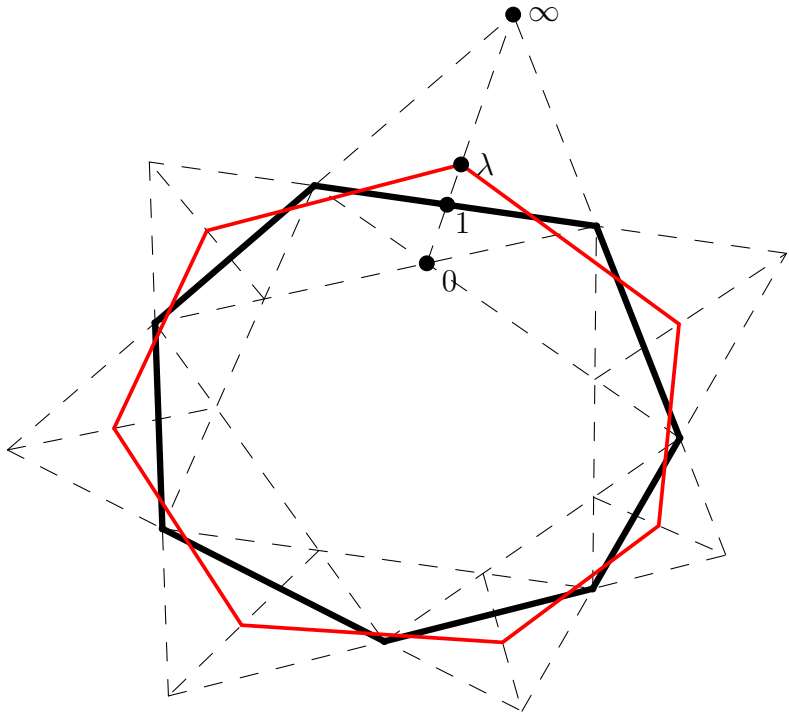
\begin{figure}[!ht]
	\centering
	\begin{tikzpicture}[line cap=round,line join=round,>=triangle 45,x=1.0cm,y=1.0cm,scale = 0.8]
	\draw [line width=2.4pt] (-2.6021129085678947,0.5989637701741396)-- (0.046226113444625745,2.880609696831082);
	\draw [line width=2.4pt] (0.046226113444625745,2.880609696831082)-- (4.7317489985437,2.2083390220125185);
	\draw [line width=2.4pt] (4.7317489985437,2.2083390220125185)-- (6.117034025442558,-1.315989061127223);
	\draw [line width=2.4pt] (6.117034025442558,-1.315989061127223)-- (4.670633482651104,-3.8217252127236865);
	\draw [line width=2.4pt] (4.670633482651104,-3.8217252127236865)-- (1.2074209154039615,-4.697714273850905);
	\draw [line width=2.4pt] (1.2074209154039615,-4.697714273850905)-- (-2.4798818767827013,-2.823505119811274);
	\draw [line width=2.4pt] (-2.6021129085678947,0.5989637701741396)-- (-2.4798818767827013,-2.823505119811274);
	\draw [dash pattern=on 7pt off 7pt] (-2.4798818767827013,-2.823505119811274)-- (4.670633482651104,-3.8217252127236865);
	\draw [dash pattern=on 7pt off 7pt] (4.670633482651104,-3.8217252127236865)-- (4.7317489985437,2.2083390220125185);
	\draw [dash pattern=on 7pt off 7pt] (4.7317489985437,2.2083390220125185)-- (-2.6021129085678947,0.5989637701741396);
	\draw [dash pattern=on 7pt off 7pt] (-2.6021129085678947,0.5989637701741396)-- (1.2074209154039615,-4.697714273850905);
	\draw [dash pattern=on 7pt off 7pt] (1.2074209154039615,-4.697714273850905)-- (6.117034025442558,-1.315989061127223);
	\draw [dash pattern=on 7pt off 7pt] (6.117034025442558,-1.315989061127223)-- (0.046226113444625745,2.880609696831082);
	\draw [dash pattern=on 7pt off 7pt] (0.046226113444625745,2.880609696831082)-- (-2.4798818767827013,-2.823505119811274);
	\draw [dash pattern=on 7pt off 7pt] (-5.055228130744793,-1.5144893443167295)-- (3.4940736662314396,-5.859990810183104);
	\draw [dash pattern=on 7pt off 7pt] (3.4940736662314396,-5.859990810183104)-- (7.889846181602513,1.7552207305301668);
	\draw [dash pattern=on 7pt off 7pt] (7.889846181602513,1.7552207305301668)-- (-2.697660551904481,3.2742977835985627);
	\draw [dash pattern=on 7pt off 7pt] (6.882081180923875,-3.262359030219397)-- (3.3490407324734592,5.726111522455926);
	\draw [dash pattern=on 7pt off 7pt] (-5.055228130744793,-1.5144893443167295)-- (3.3490407324734592,5.726111522455926);
	\draw [dash pattern=on 7pt off 7pt] (-5.055228130744793,-1.5144893443167295)-- (-1.588398659806022,-0.8104785008316736);
	\draw [dash pattern=on 7pt off 7pt] (-0.7879080466618661,0.9970809482035184)-- (-2.697660551904481,3.2742977835985627);
	\draw [dash pattern=on 7pt off 7pt] (1.9134118185022277,1.589870585281195)-- (3.3490407324734592,5.726111522455926);
	\draw [dash pattern=on 7pt off 7pt] (4.70591598614076,-0.3405182017441019)-- (7.889846181602513,1.7552207305301668);
	\draw [dash pattern=on 7pt off 7pt] (4.686039617314619,-2.3016532592568413)-- (6.882081180923875,-3.262359030219397);
	\draw [dash pattern=on 7pt off 7pt] (2.848489311785467,-3.56735180995327)-- (3.4940736662314396,-5.859990810183104);
	\draw [dash pattern=on 7pt off 7pt] (0.12052723160439807,-3.1865251947712965)-- (-2.380533641044462,-5.605255720481978);
	\draw [dash pattern=on 7pt off 7pt] (-2.697660551904481,3.2742977835985627)-- (-2.380533641044462,-5.605255720481978);
	\draw [dash pattern=on 7pt off 7pt] (-2.380533641044462,-5.605255720481978)-- (6.882081180923875,-3.262359030219397);
	\draw [red, line width=1.4pt] (2.4839546833847708,3.233681692158508)-- (6.104429846993197,0.5800170379075922);
	\draw [red, line width=1.4pt] (6.104429846993197,0.5800170379075922)-- (5.766872633738004,-2.7744869858727594);
	\draw [red, line width=1.4pt] (5.766872633738004,-2.7744869858727594)-- (3.171459219180338,-4.714302476422784);
	\draw [red, line width=1.4pt] (3.171459219180338,-4.714302476422784)-- (-1.163351553939196,-4.428141040217626);
	\draw [red, line width=1.4pt] (-1.163351553939196,-4.428141040217626)-- (-3.28373957390498,-1.1547522520889972);
	\draw [red, line width=1.4pt] (-3.28373957390498,-1.1547522520889972)-- (-1.740994969797888,2.1335557432369647);
	\draw [red, line width=1.4pt] (-1.740994969797888,2.1335557432369647)-- (2.4839546833847708,3.233681692158508);
	\begin{scriptsize}
	\draw [fill=black] (-2.6021129085678947,0.5989637701741396) circle (1.5pt);
	\draw [fill=black] (0.046226113444625745,2.880609696831082) circle (1.5pt);
	\draw [fill=black] (4.7317489985437,2.2083390220125185) circle (1.5pt);
	\draw [fill=black] (6.117034025442558,-1.315989061127223) circle (1.5pt);
	\draw [fill=black] (4.670633482651104,-3.8217252127236865) circle (1.5pt);
	\draw [fill=black] (1.2074209154039615,-4.697714273850905) circle (1.5pt);
	\draw [fill=black] (-2.4798818767827013,-2.823505119811274) circle (1.5pt);
	\draw [fill=black] (2.2515834829052217,2.5641888568649964) circle (3.5pt) node[below right] {\Large $1$};
	\draw [fill=black] (1.9134118185022277,1.589870585281195) circle (3.5pt) node[below right] {\ \Large $0$};
	\draw [fill=black] (3.349040732473459,5.726111522455926) circle (3.5pt) node[right] {\ \Large $\infty$};
	\draw [fill=black] (2.4839546833847708,3.233681692158508) circle (3.5pt) node[right] {\ \Large $\lambda$};
	\end{scriptsize}
	\end{tikzpicture}
	\caption{A parameterized heat map.} \label{fig:heat1P}
\end{figure}

It is obvious that this construction commutes with any projective transformation. Therefore, we can consider $\Hlam$ 
both as a map on the space $\mathcal{P}_n = (\RR P^2)^n$ of $n$-gons and as a map on the space $\mathfrak{P}_n = 
\mathcal{P}_n/PGL(3,\RR)$ of projective equivalence classes of $n$-gons. Note that $\Hlam$ is not defined everywhere,
but is defined at least on the full measure subset of generic $n$-gons.

In computer simulation, the maps $H_{\lambda}$ exhibit a similar convergence behavior as the projective heat map and 
the midpoint map. However, the limit $n$-gon might be regular or star-regular depending on $\lambda$ and $n$ in a very 
concrete way. We suggest the following conjecture on the dynamics of  $H_{\lambda}$ as maps on $\mathfrak{P}_n$.

\begin{conj}\label{conj:regular}	
	Let $n \geq 5$. If $\lambda \in (0, \infty)$, then $\lim_{k \to \infty} H_{\lambda}^k (P) = P_{reg}$, for almost all
	$P$ in $\mathfrak{P}_n$. Here $P_{reg}$ is the equivalence class of the regular polygon.
	
	Suppose $\lambda \in (-\infty, 0)$. 
	\begin{itemize}
	  \item[(i)] If $n = 3m \pm 1$, then $\lim_{k \to \infty} H_{\lambda}^k (P) = 
		P_{m-reg}$ for almost all $P$ in $\mathfrak{P}_n$, where $P_{m-reg}$ is the equivalence class of the $m$-regular 
		polygon whose vertices are the points $\exp(2jm\pi i/n) \in \CC \subset \RR P^2$, for $j = 0,\ldots,n-1$. 
	  \item[(ii)] If $n = 3m$, then, for almost all $P$ in $\mathfrak{P}_n$, the iterates $H_{\lambda}^k (P)$ will
	  degenerate in the sense that they will escape all compact subsets of the set of generic $n$-gons.
	\end{itemize}
\end{conj}

The following theorem is the the first step toward resolving Conjecture \ref{conj:regular} for positive $\lambda$.

\begin{thm}[Preservation of convexity] \label{thm:convex}
	Let $\lambda > 0$. Then $\Hlam$ carries a convex $n$-gon to a convex one. In other words, the set $\mathcal{C}_n$ of
	convex $n$-gons  is preserved by $\Hlam$, so is the corresponding set $\mathfrak{C}_n$ in $\mathfrak{P}_n$. 
\end{thm}
Notice that $\Hlam$ is defined everywhere in $\mathcal{C}_n$ and $\mathfrak{C}_n$. Therefore, this theorem implies that, if $P$ is convex, $H^k(P)$ is always defined for any $k$.

The rest of our results only work for pentagons, due to the non-linear nature of $\Hlam$. In the case of \textit{convex} pentagons, Conjecture 1 is completely resolved.
\begin{thm}\label{thm:convergence}
	Let $P$ be an equivalence class of convex pentagons in $\mathfrak{C}_5$. Then, the iterations $\{ \Hlam^k(P)\}_{k=1}^n$ converge to the regular class if $\lambda > 0$, and to the star-regular class if $\lambda < 0$. The rates of convergence for both cases are exponential.
\end{thm}

There is also an interesting phenomenon occurring for pentagons; that is the map $\Hlam$ reduces to a constant 
function for certain parameters.
\begin{thm}\label{thm:constant}
	Consider the action of $\Hlam$ on $\mathfrak{P}_5$. Let $\phi = \frac{1}{2} (1 + \sqrt{5})$ be the golden ratio. Then, for any generic $P \in \mathfrak{P}_5$, 
	$H_{\phi}(P) = P_{reg}$ and $H_{-1/\phi}(P) = P_{\star-reg}$, where $P_{reg}$ and $P_{\star-reg}$ are the regular and star-regular class, respectively.
\end{thm}
This theorem reminds us of the fact that the midpoint map carries any quadrilateral to a affinely regular one.

If we consider the action of $\Hlam$ on the space $\mathcal{P}_n$ instead of $\mathfrak{P}_n$, other phenomena emerge. 
Guided by computer experiments, we propose another conjecture.

\begin{conj}\label{conj:collapse}
	Let $n \geq 5$ and $\lambda \in (0,\infty)$. Set $\lambda_0 = 2 \cos(\pi/n)$.
	\begin{itemize}
		\item[(i)] When $\lambda \in (0, \lambda_0)$, any $n$-gon will collapse to a point under the iteration of $\Hlam$.
		\item[(ii)] When $\lambda = \lambda_0$, for any  $n$-gon $P$, its even iterates $\Hlam^{2k}(P)$ 
		will converge (without rescaling or any other normalization) to a regular polygon and its odd iterates $\Hlam^{2k+1}(P)$ will converge to another regular polygon. The vertices of these two regular polygons form a regular $2n$-gon.
		\item[(iii)] When $\lambda \in (\lambda_0, \infty)$, under iteration of $\Hlam$, any convex $n$-gon will degenerate so
		that its vertices will all approach a straight line.
	\end{itemize} 
\end{conj}

The parameter $\lambda_0$ is the only one in $(0,\infty$) such that $\Hlam$ maps a regular $n$-gon to another one with 
the same circumcircle. As for Conjecture \ref{conj:regular}, our results are limited to the case of pentagons.
\begin{thm} \label{thm:collapse}
	The following holds.
	\begin{itemize}
		\item [(i)] Conjecture \ref{conj:collapse} is true for \textit{convex} pentagons. 
		\item [(ii)] Moreover, the point of collapse, in the case
		of $\lambda \in (0, \lambda_0)$, and the line of collapse, in the case of $\lambda \in (\lambda_0, \infty)$ vary analytically 
	according to the vertices of the original pentagon.
	\end{itemize}
\end{thm}

We expect the point of collapse to be a complicated function of the original vertices because it is a product of a
limiting  process. However, for the parameter $1/\phi$, the point of collapse can be computed explicitly. First, we need the following definition. 
\begin{defi} \label{def:center}	
	Given a generic pentagon $P$ in $\mathcal{P}_5$, we know, by Theorem \ref{thm:constant}, that $H_{\phi}(P)$ is projectively regular. Define $Center(P)$ to be the (projective) center of $H_{\phi}(P)$, i.e. the image of the origin under the projective transformation carrying the regular pentagon (inscribed in the unit circle) to $H_{\phi}(P)$.
\end{defi}

We have the following curious result.
\begin{thm} \label{thm:center}
	Given a generic pentagon $P$. Then, the point of collapse of $P$ under iteration of $H_{1/\phi}$ is precisely 
	$Center(P)$.
\end{thm}

This article is organized as following. Sections 2 and 3 describe the basic definitions and the method of positive dominance, our main tool. The next three sections include the proofs of our results. The last section lists some ideas for further investigation.

{\bf Acknoledgement.} The author would like to express his deepest gratitude to his advisor, Professor Richard E. Schwartz, for his guidance and support throughout the whole project. He would also like to thank Karim Adiprasito, whose grant \textit{ERC StG 716424 - CASe} supports part of this work.

%

\section{Basic definitions}
\subsection{Definition of the heat maps} \label{sec:basic}
\begin{figure}[!ht]
	\centering
	\begin{tikzpicture}[line cap=round,line join=round,>=triangle 45,x=1.0cm,y=1.0cm]
	\draw [line width=2.4pt] (-1.54,-0.66)-- (0.06,1.98);
	\draw [line width=2.4pt] (0.06,1.98)-- (3.26,2.38);
	\draw [line width=2.4pt] (3.26,2.38)-- (5.92,0.24);
	\draw [dash pattern=on 5pt off 5pt] (-1.54,-0.66)-- (3.26,2.38);
	\draw [dash pattern=on 5pt off 5pt] (5.92,0.24)-- (0.06,1.98);
	\draw (1.8086119603766662,1.4607875749052217)-- (1.2718241690917447,3.979509879001378);
	\draw [dash pattern=on 5pt off 5pt] (0.06,1.98)-- (1.2718241690917447,3.979509879001378);
	\draw [dash pattern=on 5pt off 5pt] (3.26,2.38)-- (1.2718241690917447,3.979509879001378);
	\begin{scriptsize}
	\draw [fill=black] (-1.54,-0.66) circle (1.5pt) node[left] {\Large $A_0$\ \ };
	\draw [fill=black] (0.06,1.98) circle (1.5pt) node[left] {\Large $A_2$\ \ };
	\draw [fill=black] (3.26,2.38) circle (1.5pt) node[right] {\Large\ \ $A_4$};
	\draw [fill=black] (5.92,0.24) circle (1.5pt) node[right] {\Large\ \ $A_6$};
	\draw [fill=black] (1.8086119603766662,1.4607875749052217) circle (3pt) node[below right] {\large\ $0$}  node[below left] {\Large $P$\ \ };
	\draw [fill=black] (1.2718241690917447,3.979509879001378) circle (3pt) node[above right] {\large\ $\infty$}  node[left] {\Large $Q$\ \ };
	\draw [fill=black] (1.6554551413713239,2.179431892671415) circle (3pt) node[below right] {\large\ $1$} node[below left] {\Large $H$\ };
	\draw [fill=red] (1.5247990525608937,2.792498006312071) circle (3pt) node[right] {\large\ $\lambda$} node[left] {\Large $B_3$\ \ };
	\end{scriptsize}
	\end{tikzpicture}
	\caption{The heat map $H_{\lambda}$.} \label{fig:heat1P-def}
\end{figure}
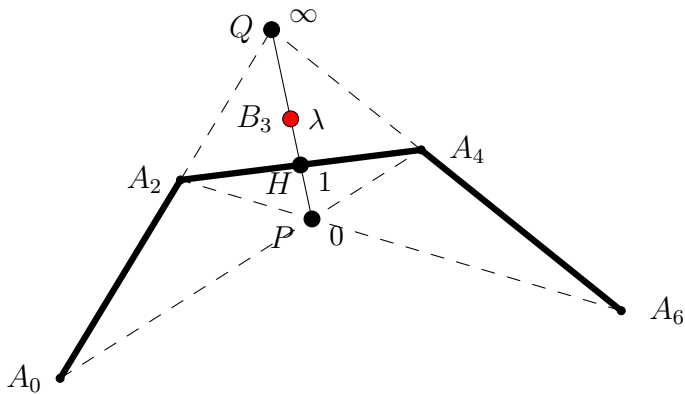

Let $\lambda \in \RR \cup {\infty}$. Starting with an $n$-gon $A$ with vertices $A_0, A_2, ..., A_{2n}$, we
construct an $n$-gon $B = \Hlam(A)$ whose vertices are labeled as $B_1, B_3, ..., B_{2n+1}$. The alternating 
labeling will be useful for later calculations. Figure \ref{fig:heat1P-def} shows the construction.

Here, the collinear points $P, Q, H$ are those that appear in the construction of the pentagram map, the inverse of the
pentagram map and Schwartz' projective heat map, respectively. Also, $B_3$ is positioned so that the cross-ratio $[H,
Q, P, B_3] = \lambda$, where, following \cite{S1}, we also define \[ [a,b,c,d] = \frac{(a-b)(c-d)}{(a-c)(b-d)}.
\] When $a < b < c <d$, this quantity lies in $(0,1)$. This occurs frequently when we consider convex polygons.

Another way to describe the above construction is as following. Consider $\RR^2$ as the usual affine patch of $\RR P^2$ 
with its usual inhomogeneous coordinates. There exists a unique projective transformation $T \in PGL(3, \RR)$ such that 

$$T(A_0) = (-1, 1),\qquad T(A_2) = (1,1), $$ $$ T(A_4) = (1,-1), \quad \text{and} \quad T(A_6) = (-1,-1).$$ Then, $T(P)$ is the origin and $T(H) =
(1,0)$, while $T(Q)$ is the point at infinity of the $x$-axis. We take $T(B_3) = (\lambda,0)$. It is easily seen that
$[H, Q, P, B_3] = [T(H), T(Q), T(P), T(B_3)] = \lambda$. Therefore, we can define $B_3 := T^{-1} (\lambda, 0)$.

We repeat the process to construct the rest of the points of $B$. This yields an $n$-gon $B = \Hlam(A)$. The
construction is projectively natural and can be carried out over any field with characteristic different from $2$. 
However, we will be primarily interested in $\RR P^2$.

\begin{figure}
	\centering
	\begin{tikzpicture}[line cap=round,line join=round,>=triangle 45,x=1.0cm,y=1.0cm, scale=1.1]
	\draw [dash pattern=on 5pt off 5pt] (-2.42,0.46)-- (1.76,5.18);
	\draw [dash pattern=on 5pt off 5pt] (1.76,5.18)-- (2.62,-0.8);
	\draw [dash pattern=on 5pt off 5pt] (2.62,-0.8)-- (-2.32,3.8);
	\draw [dash pattern=on 5pt off 5pt] (-2.32,3.8)-- (4.54,3.7);
	\draw [dash pattern=on 5pt off 5pt] (4.54,3.7)-- (-2.42,0.46);
	\draw [dash pattern=on 5pt off 5pt] (-0.7537299035369772,2.341529869690303)-- (1.9674497273846512,3.7375007328369576);
	\draw [dash pattern=on 5pt off 5pt] (1.9674497273846512,3.7375007328369576)-- (0.038035883852266444,1.604258083862262);
	\draw [dash pattern=on 5pt off 5pt] (0.038035883852266444,1.604258083862262)-- (0.5014576854121695,3.758870879221397);
	\draw [dash pattern=on 5pt off 5pt] (0.5014576854121695,3.758870879221397)-- (2.1339231475976868,2.579929741123061);
	\draw [dash pattern=on 5pt off 5pt] (2.1339231475976868,2.579929741123061)-- (-0.7537299035369772,2.341529869690303);
	\draw [line width=2.4pt] (-0.7537299035369772,2.341529869690303)-- (0.5014576854121695,3.758870879221397);
	\draw [line width=2.4pt] (0.5014576854121695,3.758870879221397)-- (1.9674497273846512,3.7375007328369576);
	\draw [line width=2.4pt] (1.9674497273846512,3.7375007328369576)-- (2.1339231475976868,2.579929741123061);
	\draw [line width=2.4pt] (2.1339231475976868,2.579929741123061)-- (0.038035883852266444,1.604258083862262);
	\draw [line width=2.4pt] (0.038035883852266444,1.604258083862262)-- (-0.7537299035369772,2.341529869690303);
	\draw [dash pattern=on 5pt off 5pt] (0.3144748459912065,2.889521148645735)-- (-2.32,3.8);
	\draw [dash pattern=on 5pt off 5pt] (1.127619907523292,3.306666272220665)-- (1.76,5.18);
	\draw [dash pattern=on 5pt off 5pt] (1.3999199157717892,3.110015446175693)-- (4.54,3.7);
	\draw [dash pattern=on 5pt off 5pt] (0.822564093837536,2.471666101315031)-- (2.62,-0.8);
	\draw [dash pattern=on 5pt off 5pt] (0.21379114537740976,2.4214068111320883)-- (-2.42,0.46);
	\draw [red, line width=1.4pt] (-0.950619641494483,1.5542601308288044)-- (-0.7776001219220474,3.2669440585252563);
	\draw [red, line width=1.4pt] (-0.7776001219220474,3.2669440585252563)-- (1.4161212363058735,4.161309351829423);
	\draw [red, line width=1.4pt] (1.4161212363058735,4.161309351829423)-- (2.836263335759666,3.379887681612106);
	\draw [red, line width=1.4pt] (2.836263335759666,3.379887681612106)-- (1.478728394912765,1.2773256002930942);
	\draw [red, line width=1.4pt] (1.478728394912765,1.2773256002930942)-- (-0.950619641494483,1.5542601308288044);
	\begin{scriptsize}
	\draw [fill=black] (-0.7537299035369772,2.341529869690303) circle (1.5pt) node[left] {\Large $A_0$\ \ };
	\draw [fill=black] (1.9674497273846512,3.7375007328369576) circle (1.5pt);
	\draw [fill=black] (0.038035883852266444,1.604258083862262) circle (1.5pt);
	\draw [fill=black] (0.5014576854121695,3.758870879221397) circle (1.5pt);
	\draw [fill=black] (2.1339231475976868,2.579929741123061) circle (1.5pt);
	\draw [fill=black] (-0.950619641494483,1.5542601308288044) circle (1.5pt);
	\draw [fill=black] (-0.7776001219220474,3.2669440585252563) circle (1.5pt);
	\draw [fill=black] (1.4161212363058735,4.161309351829423) circle (1.5pt);
	\draw [fill=black] (2.836263335759666,3.379887681612106) circle (1.5pt) node[below right] {\Large $B_0$};
	\draw [fill=black] (1.478728394912765,1.2773256002930942) circle (1.5pt);
	\end{scriptsize}
	\end{tikzpicture}
	\caption{Labeling for pentagons.} \label{fig:label-penta}
\end{figure}
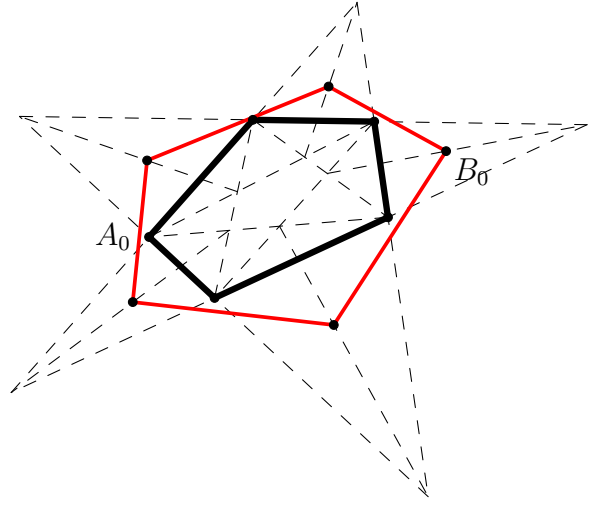

When dealing with pentagons, we prefer the labeling convention as shown in Figure \ref{fig:label-penta}. Here, $A = (A_0, A_2, A_4, A_6, A_8)$ is mapped to $B = (B_0, B_2, B_4, B_6, B_8)$.

\subsection{Flag invariants}
On an oriented polygon, a \textit{flag} is a pair $(v,e)$, where $v$ is a vertex and $e$ is an edge adjacent to $v$. We
will present the flag as an arrow which originates from $v$ and points to the next vertex on $e$. We also say that a (finite)) sequence of flags $f_n = (v_n, e_n)$ is
\textit{linked} if $e_n$ is the edge connecting $v_n$ and $v_{n+1}$. A sequence of flags is \textit{generic} if
the points $v_n$ are in general positions and so are the (unoriented) lines $e_n$. 

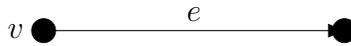
\begin{figure}[!ht]
	\centering
	\begin{tikzpicture}[line cap=round,line join=round,>=triangle 45,x=1.0cm,y=1.0cm]
	\draw [->] (0,0) --node[above] {\Large $e$} (4,0) ;
	\draw [fill=black] (0,0) circle (4.5pt) node[left] {\Large $v$\ \ };
	\draw [fill=black] (4,0) circle (4.5pt);
	\end{tikzpicture}
	\caption{The flag $(v,e)$.} \label{fig:flag1}
\end{figure}

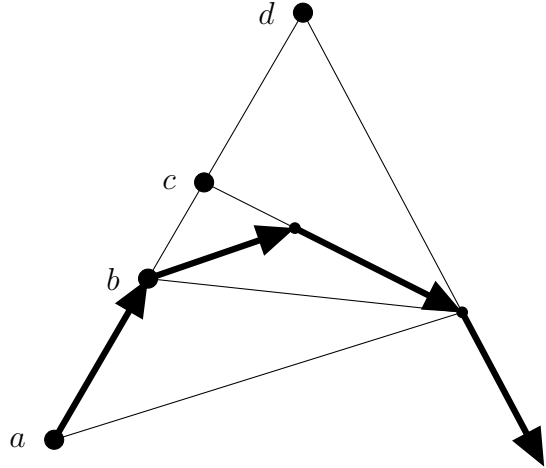
\begin{figure}[!ht]
	\centering
	\begin{tikzpicture}[line cap=round,line join=round,>=triangle 45,x=1.0cm,y=1.0cm, scale=0.8]
	\draw [->,line width=2.4pt] (-2.5,-0.08) -- (-0.94,2.6);
	\draw [->,line width=2.4pt] (-0.94,2.6) -- (1.5,3.44);
	\draw [->,line width=2.4pt] (1.5,3.44) -- (4.28,2.04);
	\draw [->,line width=2.4pt] (4.28,2.04) -- (5.642944387417847,-0.5271965013876505);
	\draw (-0.94,2.6)-- (1.6345323248829842,7.022914506850255);
	\draw (1.6345323248829842,7.022914506850255)-- (4.28,2.04);
	\draw (-0.00876691854189161,4.1998106784023905)-- (1.5,3.44);
	\draw (-0.94,2.6) -- (4.28,2.04);
	\draw (-2.5,-0.08) -- (4.28,2.04);
	\begin{scriptsize}
	\draw [fill=black] (-2.5,-0.08) circle (4.5pt) node[left] {\Large $a$\ \ \ };
	\draw [fill=black] (-0.94,2.6) circle (4.5pt) node[left] {\Large $b$\ \ \ };
	\draw [fill=black] (1.5,3.44) circle (2.5pt);
	\draw [fill=black] (4.28,2.04) circle (2.5pt);
	\draw [fill=black] (-0.00876691854189161,4.1998106784023905) circle (4.5pt) node[left] {\Large $c$\ \ \ };
	\draw [fill=black] (1.6345323248829842,7.022914506850255) circle (4.5pt) node[left] {\Large $d$\ \ \ };
	\end{scriptsize}
	\end{tikzpicture}
	\caption{A quadruple of linked flags.} \label{fig:flag2}
\end{figure}

There is only one generic triple of linked flags up to projective transformations. On the other hand, the generic linked quadruple of flags $(f_1,f_2,f_3,f_4)$ can be parameterized by $\RR P^1$. The parameterization is given by the cross-ratio $[a,b,c,d]$ as
shown in Figure \ref{fig:flag2}. We call this cross-ratio the \textit{flag invariant} associated to linked quadruple and
denote it as $[f_1, f_2, f_3, f_4]$.

A generic $n$-gon $A$ give rise to $2n$ flags. They fit in two linked sequence of flags each of which has length $n$. We
indicate each such flag, say $(v,e)$, as an arrow from $v$ to the other vertex of $A$ that lies on $e$. The arrow is
only to distinguish between the $2$ flags containing the vertex $v$ and should not be thought of as an orientation of $e$.

We label the edges by $e_{2i+1} = A_{2i} A_{2i+2}$ and the flags by $f_{2i} = (A_{2i}, e_{2i+1})$ and $f_{2i+1} =
(A_{2i+2}, e_{2i+1})$, for $i = 0, \ldots, n-1$. Note that the odd-numbered flags are linked, so are the even-numbered
ones. We associate to flag $f_{i}$ the flag invariant $x_i := [f_{i-3}, f_{i-1}, f_{i+1}, f_{i+3}]$. Here, the indices
are taken module $2n$. Furthermore, if we assume the convexity of $A$, then it is clear from the definition that $0 <
x_i < 1$ for all $i = 0,\ldots,2n-1$. These are exactly the flag invariants defined in \cite{S1}, but we feel our
interpretation has a certain geometric appeal.

Notice that, in Figure \ref{fig:flag2}, the flag invariant can also be defined as the cross-ratio of four concurrent lines.
Hence, its definition is invariant not only under projective transformations, but also under projective dualities. In
particular, for a pentagon, its flag invariants are $5$-periodic, i.e. they have the form $(x_0, \ldots, x_4,
x_0, \ldots, x_4)$. This can also easily seen geometrically by tracing through our construction for pentagons. This shows that pentagons are self-dual.

\subsection{The action of the Gauss group on convex pentagons}
A generic $n$-gon can be determined by $2n-8$ consecutive flag invariants. Therefore, for pentagons, we only need $2$
consecutive invariants. We collect below the results of Section 3 of \cite{S1} in a slightly different wording.

	The space $\mathfrak{P}_5$ of projective equivalence classes of pentagons can be identified with a subset of $\RR^2$
	via the map $\Psi(P) = (x,y)$, where $(x,y) = (x_3(P), x_4(P))$. We can certainly use any other pair of consecutive flag
	invariants instead of $(x_3, x_4)$.

\begin{lem}	
	The space $\mathfrak{C}_5$ of convex classes are carried by $\Psi$ homeomorphically onto $(0,1)^2$. Also, the regular
	class corresponds to the point $(1/\phi,1/\phi)$ and the star-regular one is mapped to $(-\phi, -\phi)$, with $\phi =
	\frac{1}{2} (1+\sqrt{5})$  being the golden ratio.
\end{lem}

When $n > 5$, we get an embedding of the space of convex classes of $n$-gons into $(0,1)^{2n-8}$, which is generally not
surjective.

	For an arbitrary pentagon, its flag invariants satisfy 
	\[ x_{k+2} = \frac{1-x_k}{1-x_k x_{k+1}}. \]
	In order words,
	\[ (x_{k+1}, x_{k+2}) = G(x_k, x_{k+1}), \quad \text{with} \quad G(x,y) = (y, \frac{1-x}{1-xy}). \]
	
	The map $G$ is called the \textit{Gauss recurrence}. It is an order $5$ birational map, which, together with the
	reflection $R(x,y) = (y,x)$, generates the \textit{Gauss group} $\Gamma$. This is a dihedral group of order $10$, which
	encodes how the invariants change under a dihedral relabeling of the pentagon.

\begin{lem} \label{lem:FunDom}
	The triangle $T$ with vertices $(0,0), (1/\phi,0)$ and $(1/\phi,1/\phi)$ contains a fundamental domain for the
	action of $\Gamma$ on $(0,1)^2$.
\end{lem}

	This result can be derived directly from the explicit formula of $G$, as shown in \cite{S1}. It will be helpful
	in proving Theorem \ref{thm:convergence}.

\section{The method of positive dominance}
Our main tool to tackle the non-linearity of the heat maps is the method of positive dominance which was invented by R.
Schwartz in \cite{S1}. In that monograph, the method was used in various geometric applications. Here, we will
only use it to prove the positivity of certain polynomials. The results in this section are taken from the above-cited
monograph.
\subsection{The single-variable case}
Consider $P(x) = a_0 + a_1 x + \ldots + a_n x^n$ as a polynomial with real coefficients on $[0,1]$. Consider the partial
sums $A_k = a_0 + \ldots + a_k$. We call $P$ \textit{positive dominant} if $A_k > 0$ for all $k$.

\begin{lem}\label{lem:PosDom}
	If $P$ is positive dominant, then $P>0$ on $[0,1]$.
\end{lem}
We reproduce the proof of this simple lemma here to illustrate the basic idea behind the whole machinery.
\begin{proof}
	We use the induction on the degree of $P$. The base case $\deg(P) = 0$ is obvious. The lemma is also obvious for $x=0$.
	For $x \in (0,1]$, we have
	\begin{align*}
		P(x) &= a_0 + a_1 x + \ldots + a_n x^n \\
			&\geq a_0 x + a_1 x + \ldots + a_n x^n \\
			&= x (A_1 + a_2 x + \ldots + a_n x^{n-1}) = xQ(x) > 0.
	\end{align*}
	The last inequality is due to the fact that $Q$ is a positive dominant polynomial of degree $n-1$.
\end{proof}

The converse of the above lemma is generally false. For a counter example, simply take $P(x) = x^2 - 2x + 1.1$. 

Given an interval $I \subset \RR$, we call the pair $(P,I)$ positive dominant if $P \circ A_I$ is positive dominant, 
where $A_I$ is the affine and orientation preserving map which carries $[0,1]$ to $I$. Clearly, $P > 0$ on $I$ if 
$(P,I)$ is positive dominant. Now, we can improve Lemma \ref{lem:PosDom} using the following algorithm.

\textbf{Divide-and-Conquer algorithm}
\begin{itemize}
  \item [1.] Start with a list LIST of intervals. Initially LIST consists only of $[0,1]$.
  \item [2.] Let $I$ be the first interval on LIST. We delete $I$ from $LIST$ and test whether $(P,I)$ is positive
  dominant.
  \item [3.] If $(P,I)$ is positive dominant, then we go back to Step 2 if LIST is nonempty and otherwise halt. If not,
  then we append to LIST the two intervals obtained by cutting $I$ in half, then go back to Step 2.
\end{itemize}

The algorithm always halts and furnishes a proof that $P > 0$ on $[0,1]$. In practice, we do not really need to
know that the algorithm will theoretically halt. We simply want to know that it will halt in the cases we are
dealing with. That is, we only care about a robust sufficient condition.

In higher-dimensional case, the following variants will be useful.
\begin{itemize}
  \item [$\bullet$] We call $P$ \textit{strong positive dominant} (or \textit{SPD}) if $P$ is positive dominant in the
  above sense.
  \item [$\bullet$] We call $P$ \textit{very weak positive dominant} (or \textit{VWPD}) if the partial sums $A_k$ are
  only non-negative, i.e. $A_k \geq 0$ for all $k$.
  \item [$\bullet$] We call $P$ \textit{weak positive dominant} (or \textit{WPD}) if $P$ is VWPD and $A_n > 0$.
\end{itemize}

\begin{lem}
	If $P$ is VWPD then $P \geq 0$ on $[0,1]$. If $P$ is WPD then $P > 0$ on $(0,1)$.
\end{lem}

In practice, we only use weak and strong positive dominance. Moreover, the same divide-and-conquer algorithm can be
defined for WPD polynomials. 

\subsection{The general case}
Now let $P$ be a polynomials in $x_1,\ldots,x_k$. Given a multi-index $I = (i_1,\ldots,i_k) \in (\NN \cup \{0\})^k$, we
write $x^I = x_1^{i_1} \ldots x_k^{i_k}$. Then $P$ can be written as $P = \sum a_I X^I$ with $a_I \in \RR$. 

If $I' = (i'_1, \ldots, i'_k)$, we say $I' \leq I$ if $i'_j \leq i_j$ for all $j = 1,\ldots,k$. We call $P$
\textit{very weak positive dominant} if 
\[ \sum_{I' \leq I} a_{I'} \geq 0 \text{, for all } I. \]
We call $P$ \textit{strong positive dominant} if all the above inequalities are strict. We call $P$ \textit{weak
positive dominant} if $P$ is VWPD and the total sum $\sum a_I$ is positive.

\begin{lem}
	The following is true.
	\begin{itemize}
	  \item [1.] If $P$ is VWPD, then $P \geq 0$ on $[0,1]^k$.
	  \item [2.] If $P$ is WPD, then $P > 0$ on $(0,1)^k$.
	  \item [3.] If $P$ is SPD, then $P > 0$ on $[0,1]^k$.
	\end{itemize}
	
\end{lem}

\section{Proof of Theorem \ref{thm:convex}}
First, because convex polygons are obviously generic, the maps $\Hlam$ are always defined for those polygons. However,
it is not even clear, a priori, that the image of a convex polygon under $\Hlam$ is generic.

The pentagram map $H_0$ clearly preserves convexity. The inverse of the pentagram map $H_{\infty}$ also preserves
convexity because it commutes with the duality map, which itself preserves convexity. Let $A$ be an arbitrary convex
polygon and $\lambda$ a positive number. Consider an affine patch of $\RR P^2$ such that $H_{\infty} (A)$ is convex in the usual 
sense. We want to show that $B := \Hlam(A)$ is convex.

Suppose this is not the case. Clearly, $A$ and $H_0(A)$ are both convex. The vertices of $\Hlam(A)$ lie on line segments
which are non-intersecting arcs connecting the two boundary components of the annulus bounded by $H_{\infty} (A)$ and
$H_0(A)$. Therefore, the polygon $\Hlam(A)$ cannot be self-intersecting. The only way it can be non-convex is that it
has $3$ consecutive vertices, say $B_1, B_3$ and $B_5$, such that, as we go from $B_1$ to $B_5$, we have to make a
non-convex turn at $B_3$, i.e. the angle $\angle B_1 B_3 B_5 > \pi$. Then, by the Intermediate Value Theorem (on the
angle $\angle B_1 B_3 B_5$), there exists a parameter $0 < \lambda' < \lambda$ such that the corresponding vertices
$B'_1, B'_3$ and $B'_5$ of the polygon $B' := H_{\lambda'}(A)$ are collinear. The following lemma asserts that this 
cannot happen, and hence, infer Theorem 2.

\begin{lem}
	For any $\lambda \in (0, \infty)$ and any convex polygon $A$, $\Hlam(A)$ does not have $3$ consecutive vertices which
	are collinear.
\end{lem}

\begin{proof}
	To construct $3$ consecutive vertices of $\Hlam(A)$, we need $6$ consecutive vertices of $A$. Suppose we are
	considering the vertices $\{A_0, A_2, A_4, A_6, A_8, A_{10}\}$ of $A$ and the vertices $\{ B_3, B_5, B_7 \}$ of $B :=
	\Hlam(A)$. By applying a projective transformation, we can assume $A_0 = (-1, 1), A_2 = (1,1), A_4 = (1,-1)$
	and $A_6 = (-1,-1)$. We only need four flag invariants $\{x_1, x_2, x_3, x_4\}$ to determine $A_8$ and $A_{10}$. These
	flag invariants all lie in the unit interval $(0,1)$ because of the convexity assumption. Then, we can compute the
	homogeneous coordinates of $B_3, B_5, B_7$.
	
	\[ B_3 = [\lambda : 0 : 1]\]
	\begin{align*} 
		B_5 = [& \lambda  x_1 x_2 - x_1 x_2 - \lambda x_2 + 1 \\
				&: \lambda  x_1 x_2 - x_1 x_2 + \lambda x_2 + 1 \\
				&: \lambda  x_1 x_2 - x_1 x_2 - \lambda x_2 - 1] 
	\end{align*} 
	\begin{align*} 
		B_7 = [&x_1 x_2 - \lambda x_2 x_3 x_4 + x_2 x_3 x_4-x_2 - \lambda x_4 - 1 \\
				&: x_1 x_2 + \lambda x_2 x_3 x_4 - x_2 x_3 x_4+x_2 - \lambda x_4 - 1 \\
				&: x_1 x_2 - \lambda x_2 x_3 x_4 + x_2 x_3 x_4-x_2 + \lambda x_4 + 1] 
	\end{align*}
	
	Then $B_3, B_5, B_7$ are collinear if and only if the determinant of the matrix $M$ formed by their homogeneous
	coordinates is zero. If we write $B_5 = [a:b:c]$ and $B_7 = [d:e:f]$, then $\det(M) = \lambda (bf - ce) + (ae - bd)$.
	
	\begin{align*}
		ae - bd &= 2x_2 (1 - x_1 x_2 - x_3 x_4 + x_1 x_2 x_3 x_4 + \lambda + \lambda x_3 x_4 - 2 \lambda x_1 x_2 x_3 x_4 \\ 
			&\qquad + \lambda^2 x_4 + \lambda^2 x_1 x_2 x_3 x_4) \\
			&= 2 x_2 \big( (1 - x_1 x_2) (1 - x_3 x_4) + \lambda (1 - x_1 x_2 x_3 x_4) + \lambda x_3 x_4 (1 - x_1 x_2) \\
			&\qquad + \lambda^2 x_4 + \lambda^2 x_1 x_2 x_3 x_4  \big) > 0 \\
		bf - ce &= 2 x_1 x_2 (x_2 - x_2 x_3 x_4 + \lambda - \lambda x_4 + 2 \lambda  x_2 x_3 x_4 + \lambda^2 x_4 - \lambda^2
		x_2 x_3 x_4) \\
			&= 2 x_1 x_2 \big( x_2 (1 - x_3 x_4) + \lambda (1 - x_4) + 2 \lambda x_2 x_3 x_4 + \lambda^2 x_4 (1 - x_2 x_3) \big)
			> 0
	\end{align*}
	
	The inequalities easily follow from the assumptions $0 < x_1, x_2, x_3, x_4 < 1$ and $\lambda > 0$. Hence, $\det(M) >
	0$ and we obtain our desired contradiction. Also, we would like to note that the method of positive dominance can also
	be used to establish the positivity of $\det(M)$. Despite the high number of variables, the method works pretty fast
	due to the low degrees of the variables.
\end{proof}

\section{The case of pentagons}

\subsection{Formula for the heat maps}
Under the identification of $\mathfrak{P}$ with $\RR^2$, the family of projective heat maps $\Hlam$ correspond to
following rational maps.

\[ \Hlam(x,y) = (x',y'), \]
\[ x' = \frac{\Alam(x,y) \Clam(x,y)}{\Blam(x,y) \Clam(y,x)}, \]
\[ y' = \frac{\Alam(y,x) \Clam(y,x)}{\Blam(y,x) \Clam(x,y)}, \]
\[ \Alam(x,y) = \lambda ^3 \left(x y^2-x y\right)+\lambda ^2 \left(-2 x y^2+3 x y+2 y^2-2 y-1\right) \] \[ +\lambda 
\left(2 x y^2-3 y^2+3 y-2\right)+ (y^2-y), \]
\[ \Blam(x,y) = \lambda ^3 (xy-y) + \lambda ^2 \left(2 x y^2-x y+2 x-y-2\right) \] \[ +\lambda  \left(-x y^2+4 x
y-x-2\right)+(y-1), \] 
\[ \Clam(x,y) = \lambda ^3 \left(x^2 y-x y-x+1\right) + \lambda ^2 \left(2 x^2 y^2-x^2 y-3 x
y-x+3\right) \] \[+\lambda \left(-x^2 -y^2+x^2 y-2 x y+2\right) + (-x^2 y+x).\]

The above formula can be derived as following. First, we construct the pentagon $A$ with $A_0 = (-1, 1), A_2 = (1,1),
A_4 = (1,-1), A_6 = (-1,-1)$ and $x_3 = x, x_4 = y$. Then, we construct the pentagon $B = \Hlam(A)$. Recall that we use
the labeling scheme in Figure \ref{fig:label-penta} for $B$. Finally, we compute $(x',y') = (x_3(B), x_4(B))$. We carried out this calculation
in Mathematica.

It can be computed that 
\begin{equation} \label{eq:derivative}
	d\Hlam(1/\phi,1/\phi) = \frac{(\lambda - \phi)^2}{(\lambda + \phi^3)(\lambda + \phi^{-1})} I,
\end{equation}
where $I$ is the identity matrix. It is easily seen that the fraction in the above formula is smaller than $1$, for all
$\lambda > 0$. This shows that $(1/\phi,1/\phi)$ is an attracting fixed point of $\Hlam$ with $\lambda > 0$. We note that this fixed point becomes repelling when $\lambda < 0$.

\subsection{Conjugation}
Conjecture \ref{conj:regular} also suggests a certain symmetry within the parameter space between the set $\{\lambda > 0\}$ 
and $\{\lambda < 0\}$. We are only able to pinpoint this symmetry for the case of pentagons.
\begin{prop}\label{prop:conjugacy}
	The map $H_{\lambda}$ is conjugate to  $H_{-1/\lambda}$ by the map $\star$ which carries the pentagon $A$ with vertices $(A_0,A_2,A_4,A_6,A_8)$ to the pentagon with the same vertices ordered as $(A_0,A_4,A_8,A_2,A_6)$. In other words, $H_{-1/\lambda} = \star \circ \Hlam \circ \star^{-1}$.
\end{prop}

\begin{proof}
	In Figure \ref{fig:heat1P-def}, if we interchange the middle two vertices $A_2$ and $A_4$, we will, in effect, interchange $P$ and $Q$.
	Thus, the cross-ratio $\lambda = [H,Q,P,X]$ will become $[H,P,Q,X] = -1/\lambda$. Therefore, in Figure \ref{fig:label-penta}, $H_{-1/\lambda}(\star(A))$ have the same vertices as $H_{\lambda}(A)$. In addition, if we follow the labeling scheme there, we will see that $H_{-1/\lambda}(\star(A)) = \star(H_{\lambda}(A))$.
\end{proof}

The above proposition shows that the statements for $\lambda > 0$ and $\lambda < 0$ in Theorem \ref{thm:convergence} are
equivalent to each other.

\subsection{Some interesting parameters}
Given the complicated look of the polynomials $\Alam, \Blam, \Clam$, Theorem \ref{thm:constant} comes as a totally
unexpected result.
\begin{proof}[Proof of Theorem \ref{thm:constant}]
	We follow the notation in Equation 1. When $\lambda = \phi = \frac{1}{2}(1 + \sqrt{5})$, the polynomials $\Alam, \Blam,
	\Clam$ factor as following.
	\[ \Alam(x,y) = \sqrt{5} (y + \phi) (xy + \frac{1}{\phi} y - \phi), \]
	\[ \Blam(x,y) = \phi \sqrt{5} (y + \phi) (xy + \frac{1}{\phi} x - \phi), \]
	\[ \Clam(x,y) = \phi \sqrt{5} (xy - \phi) (xy + \frac{1}{\phi} x - \phi). \]
	Hence, $\Hlam(x,y)$ reduces to $(1/\phi, 1/\phi)$, which corresponds to the regular
	pentagon. The case of $\lambda = -1/\phi$ follows from Proposition \ref{prop:conjugacy}.
\end{proof}

We would also like to note that, when $\lambda = 1/2$, $\Hlam$ simplifies significantly to
\[ H_{1/2}(x,y) = (f(x,y), f(y,x))), \text{ where } f(x,y) = \frac{(x+3) (x y^2 + xy - 2)}{(y+3) (3 x y+y-4)}. \]

\subsection{Proof of Theorem \ref{thm:convergence}}
Our proof follows the ideas of Section 4 of \cite{S1}.

For a pentagon $P = P(x,y)$, consider the quantity
\[ E(P) = E(x,y) = x_0 x_1 x_2 x_3 x_4 = \frac{xy(1-x)(1-y)}{1-xy}. \]
Here, $x = x_3, y = x_4$ and $x_0, x_1, x_2$ are the previous three flag invariants of $P$. Geometrically, when $P$ is
convex, $E(P)$ is 
a monotone function of the perimeter of $H_0(P)$ in the Hilbert metric
defined in the convex domain bounded by $P$. See Section \ref{sec:Hilbert} for a precise definition of the Hilbert metric. Note that $E(P)$ is $\Gamma$-invariant. 

For the rest of this subsection, we constrain ourselves inside the open unit square $(0,1)^2$, which corresponds to
$\mathfrak{C}_5$. We also assume $\lambda > 0$.
\begin{lem} \label{lem:Hilbert}
	If $P \in \mathfrak{C}_5$ is not the regular class, then $E(\Hlam(P)) > E(P)$ for any $\lambda > 0$.
\end{lem} 

Assuming this lemma, we can complete the proof of Theorem \ref{thm:convergence} as following. First, note that 
\[ \lim_{(x,y) \to \partial (0,1)^2} E(x,y) = 0. \]
The only issue when calculating the above limit is the point $(1,1)$ on the boundary. We resolve this issue as below.
\[ \lim_{(x,y) \to (0,0)} E(1-x,1-y) = \lim_{(x,y) \to (0,0)} \frac{(1-x)(1-y)xy}{x+y-xy} \] \[= \lim_{(x,y) \to
(0,0)} \frac{(1-x)(1-y)}{1/x + 1/y - 1} = \frac{1}{\infty} = 0. \]
Note that all the above limits are taken from within $(0,1)^2$. 

Inside $(0,1)^2$, $E(x,y)$ is obviously a continuous with values in $(0,1)$. Therefore, Lemma \ref{lem:Hilbert} implies that
$P_0 := (1/\phi, 1/\phi)$ is the unique point where $E(x,y)$ attains its maximum value. This in turn shows that the
sequence $\{ \Hlam^k(P) \}$ must converge to $P_0$ for any $P \in \mathfrak{C}_5$ and any $\lambda > 0$. When
$\{ \Hlam^k(P) \}$ gets sufficiently close to $P_0$, the attracting effect of $d\Hlam(P_0)$ forces the convergence to be exponential. Indeed, for a fixed positive $\lambda$, we can pick a constant $c_{\lambda}$ strictly between the
constant in Equation \ref{eq:derivative} and $1$. Then, for sufficiently large $k$, we have $\| \Hlam^{k+1}(P) - P_0\| <
c_{\lambda} \| \Hlam^k(P) - P_0\| $ with $\| \cdot \|$ being the Euclidean distance on $\RR^2$. This is the exponential convergence.

Now we only have to prove Lemma \ref{lem:Hilbert}. First, we write 
\[ E(\Hlam(P)) - E(P) = \frac{N(x,y;\lambda)}{D(x,y,\lambda)}, \]
where $N(x,y;\lambda), D(x,y,\lambda)$ are two reduced polynomials in the variables $x,y$ and $\lambda$. The composition
of $D(x,y;\lambda)$ is not important to us. We will focus on $N(x,y;\lambda)$. This polynomial has degree $14$ in
$\lambda$ with a zero constant coefficient, so we rewrite 
\[ N(x,y;\lambda) = \sum_{i=0}^{13} N_i(x,y) \lambda^{i+1}. \]

\begin{lem}
	$D(x,y;\lambda) > 0$ provided that $N(x,y; \lambda) > 0$ on $((0,1)^2 - (1/\phi,1/\phi)) \times (0,\infty)$.
\end{lem}

\begin{proof}
	Since $|E(\Hlam(P)) - E(P)| \leq 1$, $D(x,y; \lambda)$ can only vanish when $N(x,y; \lambda)$ does. Moreover, because
	$N$ only vanishes at $P_0$, then either $N$ and $D$ always have the same sign or always have opposite signs, inside the
	square $(0,1)^2$. We check this at $(1/2,1/2)$.
	
	Using Mathematica, we found
	\begin{align*}
		N(1/2,1/2;\lambda) &= \frac{1}{1048576} (\lambda +1) (3 \lambda +2) (4 \lambda +3)  \left(2 \lambda ^2+8 \lambda +3\right)^2 \\ &\qquad \left(4 \lambda ^6+36 \lambda ^5+91 \lambda ^4+160 \lambda ^3+169 \lambda ^2+86 \lambda +16\right), \\
		D(1/2,1/2;\lambda) &= \frac{3}{262144} (\lambda +1) (3 \lambda +2)^2 \left(\lambda ^2+6 \lambda +3\right) \\ &\qquad \left(2 \lambda ^2+8 \lambda +3\right)^2 \left(2 \lambda ^3+12 \lambda ^2+13 \lambda +4\right)^2.
	\end{align*}
	Both are clearly positive when $\lambda > 0$.
\end{proof}

Because of this lemma, Lemma \ref{lem:FunDom} and the positivity of $\lambda$, we are left to prove the following lemma.

\begin{lem}
	We have $N_i(x,y) > 0$ on the triangle $T$ of Lemma \ref{lem:FunDom}, for all $i=0,\ldots,13$.
\end{lem}
\begin{proof}
	The polynomials $N_i$ is nonnegative on the whole square $(0,1)^2$. However, thanks to Lemma \ref{lem:FunDom}, we only need
	to show the positivity of $N_i$ on $T$. Recall that $T$ is the right-angled triangle with vertices at $(0,0), (1/\phi,1/\phi)$ and $(1/\phi, 0)$.

	Some of the polynomials $N_i$ contain some of the factors $(1-x), (1-y), x, y$ and/or $(1-xy)$, and maybe some
	positive constant factors. These factors are obviously positive in the interior of $T$. Those $N_i$ are
	\begin{align*}
		N_0(x,y) &= 2 (1 - x) x (1 - y) y (1 - x y)^2 \big( x - x^2 + y - 7 x y + 6 x^2 y + x^3 y - y^2 \\
				& + 6 x y^2 - 2 x^2 y^2 - 5 x^3 y^2 - x^4 y^2 + x y^3 - 5 x^2 y^3 + 6 x^3 y^3 + x^4 y^3 \\
				& - x^2 y^4 + x^3 y^4 - x^4 y^4\big), \\
		N_1(x,y) &= (1 - x y) \big(4 x^2 + 12 x y + 4 y^2 + \ldots -x^8 y^8+7 x^8 y^7+7 x^7 y^8 \big),	\\
		N_{12}(x,y) &= (1 - x) x (1 - y) y \big( 1 + 3 x +3 y + \ldots -4 x^8 y^6-12 x^7 y^7-4 x^6 y^8 \big), \\
		N_{13}(x,y) &= 2 (1 - x)^2 x^2 (1- y )^2 y^2 (1- x y) \big( 1 - x + x^2 - y - 6 x y + 5 x^2 y - x^3 y \\
				& + y^2 + 5 x y^2 + 2 x^2 y^2 - 4 x^3 y^2 - x y^3 - 4 x^2 y^3 + x^3 y^3 + x^4 y^3 + x^3 y^4 \big),
	\end{align*}
	
	We could simply divide out those simple factors and obtain new polynomials $N_i^+$. For the other $i = 2,\ldots,11$, we
	set $N_i^+ := N_i$. The degrees in $x$ and in $y$ of these polynomials are at most $10$. On the boundary of $T$,
	most/all of the $N_i^+$ only have zeros at $(0,0)$ and/or $(1/\phi,1/\phi)$, which might greatly dampen the speed of 
	the positive dominance method. To overcome this, we subdivide the triangle $T$ into two triangles $\Delta_1, 
	\Delta_2$ and a square $\Box$.
	
\begin{figure}[!ht]
	\centering
	\begin{tikzpicture}[line cap=round,line join=round,>=triangle 45,x=1.0cm,y=1.0cm]
	\draw [line width=2.4pt] (0,0) -- (6,0);
	\draw [line width=2.4pt] (6,6) -- (6,0);
	\draw [line width=2.4pt] (6,6) -- (0,0);
	\draw (6,6) -- (0,6);
	\draw (0,0) -- (0,6);
	\draw (3,3) -- (3,0);
	\draw (3,3) -- (6,3);
	
	\draw (0,0) -- (3,2.7);
	\draw (0,0) -- (3,2.4);
	\draw (0,0) -- (3,2.1);
	\draw (0,0) -- (3,1.8);
	\draw (0,0) -- (3,1.5);
	\draw (0,0) -- (3,1.2);
	\draw (0,0) -- (3,0.9);
	\draw (0,0) -- (3,0.6);
	\draw (0,0) -- (3,0.3);
	
	\draw (6,6) -- (3.3,3);
	\draw (6,6) -- (3.6,3);
	\draw (6,6) -- (3.9,3);
	\draw (6,6) -- (4.2,3);
	\draw (6,6) -- (4.5,3);
	\draw (6,6) -- (4.8,3);
	\draw (6,6) -- (5.1,3);
	\draw (6,6) -- (5.4,3);
	\draw (6,6) -- (5.7,3);
	
	\draw (6,4.5) node[right] {\Large\ $\Delta_2$};
	\draw (1.5,0) node[below] {\Large $\Delta_1$};
	\draw (6,1.5) node[right] {\Large\ $\Box$};
	
	\draw (0,0) node[below left] {\Large $(0,0)$};
	\draw (6,6) node[above right] {\Large $\left( \frac{1}{\phi}, \frac{1}{\phi} \right)$};
	\end{tikzpicture}
	\caption{The triangle $T$.} \label{fig:FunDom}
\end{figure}
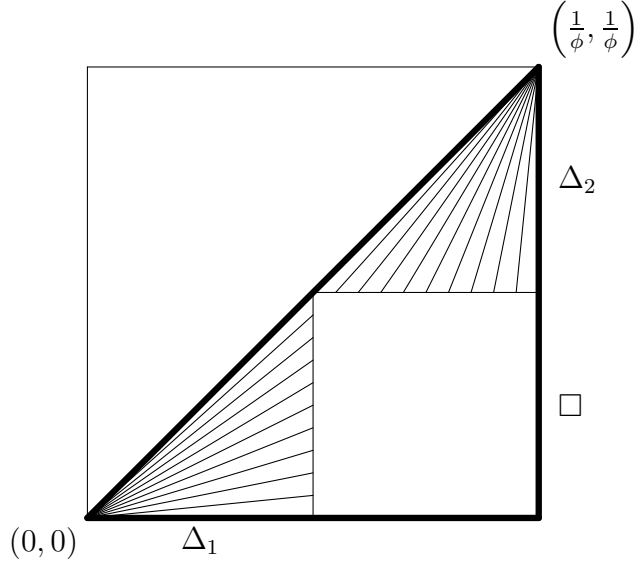
	
	For each of the triangles $\Delta_1, \Delta_2$, we map it affinely onto the triangle $T_1$ with vertices at $(0,0),(1,1)$ and $(1,0)$ so that the zero of $N_i$ is moved to
	the origin. Then we blow up the origin by the coordinate change $(x,t) \mapsto (x, tx)$ for $(x,t) \in (0,1)^2$. In effect, we
	are considering the values of $N_i^+$ along the rays through the zero. Hence, this should reveal a factor of $N_i^+$
	which is a monomial in $x$. Finally, after dividing out the monomial factor, we can apply the WPD algorithm and
	conclude our proof of this lemma. Concretely, from $N_i^+$, we create
	\[ N_i^{\Delta_1}(x,y) := N_i^+\left(\frac{x/2}{\phi}, \frac{xy/2}{\phi} \right), \] 
	\[ N_i^{\Delta_2}(x,y) := N_i^+\left(\frac{1-x/2}{\phi}, \frac{1-xy/2}{\phi} \right), \]
	\[ N_i^{\Box}(x,y) := N_i^+\left(\frac{1-x/2}{\phi}, \frac{y/2}{\phi} \right). \]  
	
	We illustrate our process for the $N_0$.
	\begin{align*}
		N_0^+(x,y) &= \big( x - x^2 + y - 7 x y + 6 x^2 y + x^3 y - y^2 + 6 x y^2 - 2 x^2 y^2 - 5 x^3 y^2 \\ &\quad - x^4 y^2 + x y^3 - 5 x^2 y^3 + 6 x^3 y^3 + x^4 y^3 - x^2 y^4 + x^3 y^4 - x^4 y^4\big) \\
		N_0^{\Delta_1}(x,y) &= x(...) \\
		N_0^{\Delta_2}(x,y) &= x^2(...) \\
		N_0^{\Delta_1}(x,y) &= (...).
	\end{align*}
	Note the factors $x$ and $x^2$ in $N_0^{\Delta_1}$ and $N_0^{\Delta_2}$. They correspond to the zeros of $N_0^+(x,y)$ at $(0,0)$ and $(1/\phi,1/\phi)$.
	
\end{proof}

	\textbf{Remark.} We note that the method of positive dominance can be applied directly to the polynomial
	$N(x,y;\lambda)$, with some tricks to handle its zeros. However, the higher dimension, together with the high degrees
	in $x,y,\lambda$,  would slow down the method significantly. The process described above is meant to alleviate these
	issues. In reality, for all $i$, the WPD algorithm runs very fast for $N_i^{\Delta_1}, N_i^{\Delta_2}, N_i^{\Box}$ and never reaches a recursive depth more than $2$.

\section{Points of collapse}
Previously, we have focused on projective equivalence classes of polygons, but for this section, we mainly consider the polygons as they are, without reducing modulo projective transformations. To fix our notation, for a polygon $P$, we will use
$[P]$ to denote its equivalence class.

\subsection{Hilbert metric} \label{sec:Hilbert}
\begin{figure}[!ht]
	\centering
	\begin{tikzpicture}[line cap=round,line join=round,>=triangle 45,x=1.0cm,y=1.0cm]
	\fill[color=black,fill=black,fill opacity=0.1] (-2.16,-0.38) -- (1.12,-1.4) -- (4.04,1.08) -- (3.68,3.74) -- (0.06,4.68) -- (-2.38,3.76) -- cycle;
	\draw [color=black] (-2.16,-0.38)-- (1.12,-1.4);
	\draw [color=black] (1.12,-1.4)-- (4.04,1.08);
	\draw [color=black] (4.04,1.08)-- (3.68,3.74);
	\draw [color=black] (3.68,3.74)-- (0.06,4.68);
	\draw [color=black] (0.06,4.68)-- (-2.38,3.76);
	\draw [color=black] (-2.38,3.76)-- (-2.16,-0.38);
	\begin{scriptsize}
	\draw [fill=black] (-1.56,1.26) circle (2.5pt);
	\draw[color=black] (-1.7,1.66) node {\Large $b$};
	\draw [fill=black] (0.56,3.38) circle (2.5pt);
	\draw[color=black] (0.32,3.72) node {\Large $c$};
	\draw [fill=black] (-2.212477064220184,0.6075229357798163) circle (2.5pt);
	\draw[color=black] (-2.46,0.72) node {\Large $a$};
	\draw [fill=black] (1.488947368421053,4.308947368421053) circle (2.5pt);
	\draw[color=black] (1.5,4.76) node {\Large $d$};
	\draw (-2.212477064220184,0.6075229357798163) -- (1.488947368421053,4.308947368421053);
	\end{scriptsize}
	\end{tikzpicture}
	\caption{The Hilbert metric.} \label{fig:Hilbert}
\end{figure}
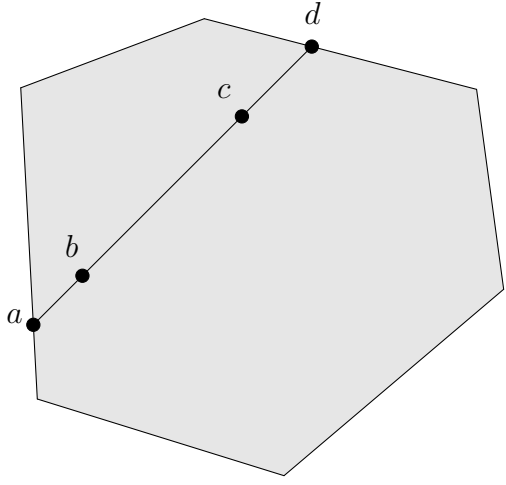

Let $K \subset \RR P^2$ be a compact convex domain. The Hilbert metric on $K$ is defined to be, for two given points
$b,c \in K$ 
\[ d_K(b,c) = -\log[a,b,c,d], \]
where $a$ and $d$ are the intersections of the line $bc$ and the boundary of $K$, ordered as in Figure \ref{fig:Hilbert}. This metric is projectively natural on $K$. When $K$ is a circle, $d_K$ is the hyperbolic metric in the Klein model. Also, if $K$ is a compact convex domain contained in the interior of another compact domain, we sometimes write $diam_{L} (K)$ for the diameter of $K$ with respect to the Hilbert metric on $L$.

We quote the two following lemmas and their proofs from \cite{S2}. 

\begin{lem}
	Suppose $\{ L_n \}$ is a sequence of compact convex subsets contained in the interior of $K$. Suppose that the Euclidean
	diameter of $L_n$ converges to the Euclidean diameter of $K$. Then, the diameter of $L_n$ defined by the Hilbert
	metric $d_K$ converges to $\infty$.
\end{lem}
\begin{proof}
	Suppose $b_n, c_n \in L_n$ are two points that realize the Euclidean diameter of $K$. Let $a_n, d_n \in \partial K$ be
	the points in the definition of $d_K(b_n,c_n)$. By our assumption on the Euclidean diameters of $L_n$ and $K$, $\|a_n
	- b_n\|$ and $\|c_n - d_n\|$ converge to $0$, whereas $\|b_n - c_n\|$ is uniformly bounded away from $0$. Therefore,
	$d_K(b_n, c_n) \to \infty$.
\end{proof}

\begin{lem}\label{lem:subset1}
	Suppose that $K_1 \supset K_2 \supset K_3 \supset \ldots $ is a nested family of compact convex subsets with $K_{n+1}$
	contained in the interior of $K_n$. Suppose that there is uniform upper bound $C$ for the diameter of $K_{n+1}$
	with respect to the Hilbert metric of $K_n$. Then $\cap K_n$ is a single point.
\end{lem}
\begin{proof}
	We just use the upper bound on diameter. The previous lemma implies that the quotient of the Euclidean diameter of
	$K_{n+1}$ and that of $K_n$ is at most $C'$, for some uniform $C' < 1$. 
\end{proof}

We have a similar result with the subsets nested in the other direction.
\begin{lem}\label{lem:subset2}
	Suppose that $K_1 \subset K_2 \subset K_3 \subset \ldots $ is a nested family of compact convex subsets with $K_{n}$
	contained in the interior of $K_{n+1}$. Suppose that there is uniform upper bound $C$ for the diameter of $K_{n}$
	with respect to the Hilbert metric of $K_{n+1}$. Then $\cup K_n$ is an affine patch of $\RR P^2$, i.e., its
	complement in $\RR P^2$ is a single line.
\end{lem}
\begin{proof}
	We associate to a compact convex subset $K$ the closure of the set of lines that are disjoint from $K$, which is
	denoted by $K^*$. Clearly, $K^*$ is a compact convex subset of the projective dual space $(\RR P^2)^*$. Note that the
	dual operator reverses the order by inclusion. We still need to show that the Hilbert diameter of $K^*_{n+1}$ with
	respect to $K^*_n$ is bounded above; which is done by the next lemma. After that, Lemma \ref{lem:subset1} shows us that
	there is exactly one line that never intersects the sets $K_n$.
\end{proof}

\begin{lem}
	If $K \subset L$ are compact convex subsets, with $K$ contained in the interior of $L$, then $diam_{K^*}(L^*) \leq
	diam_{L} (K)$.
\end{lem}
\begin{proof}
	Let $\beta, \gamma \in L^*$ be two lines where $diam_{K^*}(L^*)$ is attained. Let $\alpha, \delta$ be two lines in the
	definition of $d_{K^*}(\beta, \gamma)$, i.e. $\alpha, \delta$ pass through the intersection of $\beta$ and $\gamma$ and
	$d_{K^*}(\beta, \gamma) = -\log[\alpha, \beta, \gamma, \delta]$. Since $\alpha, \delta$ belongs to the boundary of $K^*$, it
	contains two points $a,d$ of $K$. Let $b,c$ be the intersections of $\beta, \gamma$ with the line $ad$. Take $e,f \in
	\partial L$ such that $a,d,e,f$ are collinear and $d_L{K} (a,d) = -\log[e,a,d,f]$.

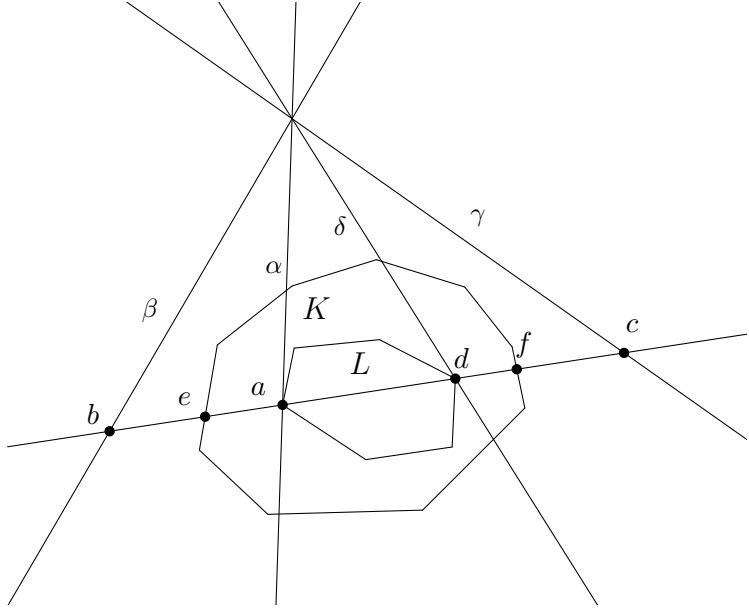
\begin{figure}[!ht]
	\centering
	
	\begin{tikzpicture}[line cap=round,line join=round,>=triangle 45,x=1.0cm,y=1.0cm, scale = 0.7]
	\clip(-4.40294,-4.71928) rectangle (9.63306,6.72732);
	\draw [domain=-4.40294:9.63306] plot(\x,{(-8.6416-5.28*\x)/-3.08});
	\draw [domain=-4.40294:9.63306] plot(\x,{(--27.3848-3.7*\x)/5.24});
	\draw (1.,1.34)-- (-0.42,0.22);
	\draw (-0.42,0.22)-- (-0.76,-1.78);
	\draw (-0.76,-1.78)-- (0.54,-3.);
	\draw (0.54,-3.)-- (3.48,-2.92);
	\draw (3.48,-2.92)-- (5.42,-0.98);
	\draw (5.42,-0.98)-- (5.18,0.18);
	\draw (5.18,0.18)-- (4.28,1.32);
	\draw (4.28,1.32)-- (2.6,1.84);
	\draw (2.6,1.84)-- (1.,1.34);
	\draw (0.82,-0.92)-- (1.04,0.16);
	\draw (1.04,0.16)-- (2.66,0.32);
	\draw (2.66,0.32)-- (4.1,-0.42);
	\draw (4.1,-0.42)-- (4.04,-1.72);
	\draw (4.04,-1.72)-- (2.4,-1.96);
	\draw (0.82,-0.92)-- (2.4,-1.96);
	\draw [domain=-4.40294:9.63306] plot(\x,{(-4.6264--5.44*\x)/0.18});
	\draw [domain=-4.40294:9.63306] plot(\x,{(--18.952-4.94*\x)/3.1});
	\draw [domain=-4.40294:9.63306] plot(\x,{(-3.4276--0.5*\x)/3.28});
	\draw (-2.03134,1.28232) node[anchor=north west] {$\beta$};
	\draw (0.31606,1.98412) node[anchor=north west] {$\alpha$};
	\draw (1.61646,2.83112) node[anchor=north west] {$\delta$};
	\draw (4.21226,2.97632) node[anchor=north west] {$\gamma$};
	\draw (1.01634,1.29592) node[anchor=north west] {\Large $K$};
	\draw (1.91326,0.25) node[anchor=north west] {\Large $L$};
	\begin{scriptsize}
	\draw [fill=black] (0.82,-0.92) circle (2.5pt);
	\draw[color=black] (0.36446,-0.60528) node {\large $a$};
	\draw [fill=black] (4.1,-0.42) circle (2.5pt);
	\draw[color=black] (4.21226,-0.02448) node {\large $d$};
	\draw [fill=black] (-0.6519449139280125,-1.1443818466353677) circle (2.5pt);
	\draw[color=black] (-1.03914,-0.79888) node {\large $e$};
	\draw [fill=black] (5.267321646962903,-0.24205462698736238) circle (2.5pt);
	\draw[color=black] (5.39806,0.21752) node {\large $f$};
	\draw [fill=black] (7.304335050149091,0.06846570886419083) circle (2.5pt);
	\draw[color=black] (7.45506,0.58052) node {\large $c$};
	\draw [fill=black] (-2.4654880089235918,-1.4208365867261572) circle (2.5pt);
	\draw[color=black] (-2.78154,-1.08928) node {\large $b$};
	\end{scriptsize}
	\end{tikzpicture}
\end{figure}
	
	Because $\beta, \gamma$ do not intersect $L$, the collinear points $b,e,a,d,f,c$ must lie in that order. This implies $[e,a,d,f] \leq [a,b,c,d] =
	[\alpha, \beta, \gamma, \delta]$. Hence, 
	\[ diam_{L} (K) \geq -\log[e,a,d,f] \geq -\log[\alpha, \beta, \gamma, \delta] = diam_{K^*}(L^*).\]
\end{proof}

\subsection{Proof of Part (i) of Theorem \ref{thm:collapse}}
	The case $\lambda = \phi = \frac{1}{2} (1 + \sqrt{5})$ is clear. Indeed, for any generic pentagon $P$, we know
	$H_{\phi}([P]) = [P_{reg}]$ by Theorem \ref{thm:constant}. Here $P_{reg}$ is the regular pentagon inscribed in the unit
	circle; its vertices are $\{\exp(2 j \pi i/5)\}_{j=0}^4 \subset \CC \subset \RR P^2$. Suppose $\Hlam(P) = T(P_{reg})$ 
	for some projective transformation $T \in PGL(3,\RR)$. Let $P'_{reg}$ be another regular pentagon inscribed in the
	unit circle with vertices $\{\exp(2 (2j+1) \pi i/10)\}_{j=0}^4$. Together, the vertices of $P_{reg}$ and $P'_{reg}$ 
	form a regular decagon. Then, we have
	\[ H_{\phi}^k(P) = \left\{ \begin{matrix} T(P_{reg}) & \text{for odd } k, \\ T(P'_{reg}) & \text{for even } k.
	\end{matrix} \right. \]
	
	Now fix $0 < \lambda < \phi$. Let $P$ be a convex pentagon. We want to appeal to Lemma \ref{lem:subset1}, but $\Hlam(P)$
	might not be contained in $P$, even when $P$ is regular. Therefore, we will consider the ellipse circumscribing $P$, 
	which always exists because we are dealing with pentagons; it will be denoted as $Ell(P)$. 
	
	First, consider a regular
	pentagon $P_{reg}$. The ellipse $Ell(\Hlam(P_{reg}))$ is properly contained in $Ell(P_{reg})$, and hence, its diameter with
	respect to the Hilbert metric on $Ell(P_{reg})$ is a finite constant. 
	
	When $P$ is close enough to being projectively regular, it must satisfy $$Ell(\Hlam(P)) \subset Ell(P).$$ In addition,
	because the Hilbert diameter of the inner ellipse with respect to outer one varies real-analytically on $P$, it is 
	bounded above by a constant when $[P]$ is in a fixed neighborhood $U$ of $[P_{reg}]$. Due to Theorem \ref{thm:convergence}, the equivalence class of $\Hlam^k(P)$ will lie in $U$ for $k \geq k_0$. For those $\Hlam^k(P)$, the Hilbert diameter of the inner ellipse can be uniformly bounded from above. Thus, Lemma \ref{lem:subset1} will force $\Hlam^k(P)$ to shrink to a point.
	
	For $\lambda > \phi$, we can use a similar argument. In this case, $Ell(\Hlam(P_{reg}))$ properly contains 
	$Ell(P_{reg})$, so we will use Lemma \ref{lem:subset2}, instead of Lemma \ref{lem:subset1}. This concludes the proof of Part (i) in all cases.
	
\subsection{Proof of Part (ii) of Theorem \ref{thm:collapse}}
	Consider a convex pentagon $P$. Let $\lambda \in (0, \phi)$. As above, we first iterate $\Hlam$ until $[P]$ gets very
	close to $[P_{reg}]$. As shown in the proof of Lemma \ref{lem:subset1}, the quotient of the Euclidean diameter of
	$Ell(\Hlam^{k+1}(P))$ and that of $Ell(\Hlam^k(P))$ is bounded above by some constant $0 < C < 1$, when $k \geq k_0$
	for some large $k_0$. Fix a vertex of $P$ and trace it as we iterate $\Hlam$, we obtain a sequence of functions
	$a_k(P)$, which obviously vary real-analytically on $P$. Now the bound on the quotient of Euclidean diameters infers
	that $\|a_k(P) - a_{k+1}(P)\| \leq d C^{k-k_0}$, where $d$ is the Euclidean diameter of $Ell(\Hlam^{k_0}(P))$. Hence,
	for $l > k \geq k_0$, we can write
	\[ \|a_l(P) - a_k(P)\| \leq d \sum_{i=k+1}^l C^{i-k_0} \leq \frac{d}{C^{k_0}} \frac{C^m}{1-C} \]
	and deduce that $a_n(P)$ are a uniformly converging family of real-analytic functions. Hence, its limit function, which
	yields the point of collapse of $P$ under the iteration of $\Hlam$, is also real-analytic. The case $\lambda \in
	(\phi, \infty)$ is resolved similarly.

\subsection{Proof of Theorem \ref{thm:center}}
	Our strategy is to normalize an arbitrary pentagon $P$ so that its first four vertices are some pre-chosen points.
	Then, we compute $Center(P)$ as a function of the last vertex, i.e., $Center(P) = f(x,y)$ with $x,y$ being the inhomogeneous
	coordinates of that vertex. Next, we normalize $H_{1/\phi}(P)$ in the same way and use the function $f$ to get its
	projective center in the normalized picture. Then we pull the point back using the normalization map and verify that it is now
	exactly $Center(P)$.
	
	Here, we opt to map the first four vertices to the points $[1:0:0], [0:1:0], [0:0:1]$ and $[1:1:1]$. The reason is the
	normalization map can be computed easily for this quadruple. Let $(x,y)$ denote the coordinates of the fifth point.
	Using Mathematica, we compute
	
	\begin{equation} \label{eq:center}
		f(x,y) = \left( \frac{\phi x + y}{-1 + \phi y}, \frac{\phi x + y}{\phi (-1 + \phi x)} \right). 
	\end{equation}
	
	The pentagon $H_{1/\phi}(x,y)$ has vertices with homogeneous coordinates
	\[ [1 : -\phi : 1], \quad [\phi^2 x: \phi x + y : \phi^3 x - \phi], \quad [\phi x + y, \phi^2 y, \phi y + 1], \]
	\[ [\phi^3 x - \phi : -\phi y + 1 : \phi^2], \quad [\phi^2 : -\phi y: 1]. \]
	
	We then compute the normalization map $T$ for $H_{1/\phi}(x,y)$ and then the projective center $c_0$ of
	$T(H_{1/\phi})(x,y)$. Finally, we verify that $T^{-1}(x_0)$ is exactly the one given in Equation \ref{eq:center}. This is our proof.
	
	We note that the entries of $T$ as a $3$-by-$3$ matrix, as well as the homogeneous coordinates of $c_0$ are rational
	functions whose degrees in $x$ and in $y$ are small (no bigger than $4$), but their coefficients might be complicated. 
	Besides, their expressions do not look enlightening, so we decide not to write them down here.
	
\subsection{Negative parameters}
In view of Proposition \ref{prop:conjugacy}, we have the following corollaries of Theorem \ref{thm:collapse} and Theorem \ref{thm:center}.

\begin{cor}
	Let $n = 5$ and $\lambda \in (-\infty,0)$. Set $\lambda_0^{\star} = -1/\phi. $.
	\begin{itemize}
		\item[(i)] When $\lambda \in (-\infty, \lambda_0^{\star})$, any pentagon will collapse to a point under the iteration of $\Hlam$.
		\item[(ii)] When $\lambda = \lambda_0^{\star}$, for any  pentagon $P$, its even iterates $\Hlam^{2k}(P)$ 
		will converge (without rescaling or any other normalization) to a star-regular pentagon and its odd iterates $\Hlam^{2k+1}(P)$ will converge to another star-regular pentagon. The vertices of these two star-regular pentagons form a regular decagon.
		\item[(iii)] When $\lambda \in (\lambda_0^{\star}, 0)$, under iteration of $\Hlam$, any convex pentagon will degenerate so that its vertices will all approach a straight line.
	\end{itemize} 
\end{cor}

Definition \ref{def:center} is modified as following.
\begin{defi}	
	Given a generic pentagon $P$ in $\mathcal{P}_5$, we know, by Theorem \ref{thm:constant}, that $H_{-1/\phi}(P)$ is projectively star-regular. Define $Center^{\star}(P)$ to be the (projective) center of $H_{-1/\phi}(P)$, i.e. the image of the origin under the projective transformation carrying the star-regular pentagon (inscribed in the unit circle) to $H_{-1/\phi}(P)$.
\end{defi}

\begin{cor}
	Given a generic pentagon $P$. Then, the point of collapse of $P$ under iteration of $H_{-\phi}$ is precisely 
	$Center^{\star}(P)$.
\end{cor}

\section{Further directions}
Conjecture \ref{conj:regular} has been resolved for $H_1$ acting on the projective equivalence classes of generic pentagons, not just the convex ones, in
\cite{S1}. In that monograph, the Julia set of $H_1$, which is defined as the classes of generic pentagons that
are never mapped by the iterations $H_1^k$ to a convex pentagon, is shown to have measure zero and connected. 

For our $1$-parameter family of heat maps, the Julia set varies wildly. It disappears when $\lambda = \phi$ and 
$\lambda = -1/\phi$; this is actually how we discovered Theorem \ref{thm:constant}. The Julia set looks messy and grainy for most
parameters, especially when we get close to $\lambda = 0$ and $\lambda = \infty$. This is expected because $H_0 =
H_{\infty}$ is the identity map on $\mathfrak{P}_5$. Curiously, for some parameters such as $\lambda = 1/2$ and 
$\lambda = 2\phi-3$, the picture of the Julia set seems to stabilize. However, for all parameters, except, of course,  for $0$ and $\infty$, the computer simulation suggests that the Julia set is always connected and of measure zero.

On another note, Theorem \ref{thm:constant}'s resemblance to the action of the midpoint map on quadrilaterals makes us
wonder if there is a higher-dimensional analogue in projective geometry. There is one in affine geometry. For an
arbitrary octahedron, the centers of mass of the facets will form a parallelogram, which is affinely regular, i.e.
affinely equivalent to a cube. Let us call this the \textit{centroid map}. Here, the octahedron needs not be convex. In
fact, we can just take a collection of $6$ points in general positions, endow it with the combinatorial (lattice) structure of an octahedron, and then, define the
facets and their centers of mass accordingly. Dually, the centroid map will take a generic collection of $8$ points with
the combinatorial structure of a cube to an affinely regular octahedron. In higher dimension, we have similar results
for hypercubes and cross-polytopes.

B. Khesin and F. Soloviev \cite{KS} have studied variants of the pentagram map on polygons with longer diagonals. We can
apply this modification to our heat maps. Here, all the other star-regular polygons appear. However, there do not seem to be
any obvious conjugations among these cases as in Proposition \ref{prop:conjugacy}. Moreover, this generalization is only available
to $n$-gons with $n \geq 6$. Here, even the convergence of the iterates of $H_1$ on convex hexagons is still a
conjecture.

Finally, we would like to note that our $1$-parameter family of heat maps is actually contained in a two-parameter
family. Fix a point $x_0 \in \RR P^2$. For any generic quadruple of points, we take the projective transformation that
carries them to four fixed points, e.g. $(-1,1), (1,1), (1,-1), (-1,-1)$ as in Section \ref{sec:basic}, and then use it to pull 
back $x_0$. Now we have a family of maps on generic quadruples parameterized by the $2$-dimensional $\RR P^2$. As
before, given an $n$-gon $P$, we can apply this map on quadruples of consecutive vertices $P$. This gives us a
two-parameter family of maps containing the one-parameter family we have been studying. Although there are maps in this
two-parameter family which behaves chaotically, there is an open set of the parameter space whose maps behave much like
the heat maps. We are planning on an article on these maps.

\end{document}